\documentclass{article}

\usepackage[english]{babel} 				
\usepackage{graphicx}					
\usepackage{vmargin}						

\usepackage{authblk}   
\usepackage{hyperref}						
\usepackage{subfig}							
\usepackage{amsmath}				    	
\usepackage{amsfonts}
\usepackage{amssymb}
\usepackage{amsthm}	
\usepackage{dsfont}  							

\usepackage{fancyhdr,fancybox} 			    
\providecommand{\keywords}[1]{\textbf{\textit{Key Words:  }} #1}
\providecommand{\msccodes}[1]{\textbf{\textit{MSC2010:  }} #1}
\usepackage{setspace}
\theoremstyle{definition}
\newtheorem{Theorem}{Theorem}[section] 
\newtheorem{Definition}[Theorem]{Definition}
\newtheorem{Assumption}[Theorem]{Assumption}
\newtheorem{Lemma}[Theorem]{Lemma}
\newtheorem{Proposition}[Theorem]{Proposition}
\newtheorem{Remark}[Theorem]{Remark}

\newcommand{\R}{\mathbb{R}}			

\newcommand{\E}{\mathbb{E}}     
\newcommand{\Var}{\mathrm{Var}} 
\newcommand{\p}{\mathbb{P}}

\author{Christof Henkel\thanks{henkel.christof@gmail.com}}

\affil{Mathematisches Institut der Ludwig-Maximilians-Universit\"at M\"unchen\\
  Theresienstrasse 39, 80333 M\"unchen, Germany\\}
\date{June 2016}
\title{An agent behavior based model for diffusion price processes with application to phase transition and oscillations}

\begin{document}

\maketitle
\begin{abstract}
We present an agent behavior based microscopic model for diffusion price processes. As such we provide a model not only containing a convenient framework for describing socio-economic behavior, but also a sophisticated link to price dynamics. We furthermore establish the circumstances under which the dynamics converge to diffusion processes in the large market limit. 
To demonstrate the applicability of a separation of behavior and price process, we show how herding behavior of market participants can lead to equilibria transition and oscillations in diffusion price processes. 
\end{abstract}

\vfill 

\msccodes{60F05, 60J60, 60K35, 91B69, 91D30}
\renewcommand{\thepage}{\arabic{page}}

\setlength{\skip\footins}{10mm} 
\thispagestyle{empty} 
\newpage
\setcounter{page}{1} 

\section{Introduction}
The foundation for modern financial modeling was set with the french mathematician Louis Bachelier suggesting a brownian motion based model for describing price fluctuation at the Paris stock exchange in 1900 (see Bachelier \cite{bach1900}). Since then many and various probabilistic models have been invented and were further developed. Most of these market models derive the dynamics of price processes from the interaction of market participants. Since in the literature not only a vast diversity of how the agent interaction is modeled but also on how the price process is derived  exists, a short overview of recent models is given in the following.\\

In F{\"o}llmer and Schweizer \cite{foe1993} as well as in Horst \cite{hor2005} stock prices are modeled in discrete time as sequence of temporary equilibria which emerge as a consequence of simultaneous matching of supply and demand of several agents. It is further shown that in a noise trader environment the resulting price process can be approximated by an Ornstein-Uhlenbeck process. Although they already capture some agent interaction and mimic effects, their model is rather simple leading to different shortcomings. Firstly, it lacks feedback effects of the price with respect to the agents behavior, which has been addressed and complementary elaborated in F{\"o}llmer, Horst and Kirman \cite{foe2005}. Secondly, simultaneous excess matching seems unrealistic in light of modern financial markets were orders arrive asynchronously in continuous time (see e.g. Bayraktar, Horst and Sircar \cite{bay2007}).

To account for this asynchronous order arrival Bayraktar, Horst and Sircar \cite{bay2006},\cite{bay2007} as well as Horst and Rothe \cite{hor2008} use the mathematical framework of queuing theory earlier examined by Mandelbaum, Pats et al. \cite{man1998a} and Mandelbaum, Massey and Reiman \cite{man1998b} for their models.

Also the model explained in Lux \cite{lux1995} and Lux \cite{lux1997} takes asynchronous order arrival into consideration by using a so called market maker who matches supply and demand and alters the price accordingly. In order to examine the connection between social economic behavior (e.g. mimic effect represented by herding behavior) and observable price process properties (e.g. volatility clustering, bubbles, crashes) the author differentiates the agents by type and assign specific characteristics.

It seems natural to characterize agents by their opinion. Opinion-based models range from binary (e.g. F{\"o}llmer \cite{foe1974}, Arthur \cite{art1994}, Orl{\'e}an \cite{orl1995}, Latan{\'e} and Nowak \cite{lat1997}, Weisbuch and Boudjema \cite{wei1999} and Sznajd-Weron and Sznajd \cite{szn2000}) to opinions from a continuous spectrum, which are used, for example, to describe  large social networks or ratings (see Duffant et al. \cite{def2000}, G{\'o}mez-Serrano, Graham and LeBoudec \cite{gom2012} or Weisbuch, Deffuant and Amblard \cite{wei2005}).\\

However the characteristic and the interaction of the agents is described in the respective model, it mostly can be classified in a wider sense as interacting objects with assigned states. Thereby the behavior of an object is modeled by a state transition to another, where the transition rule considers other objects. An early example of binary states is given by the model of Ising \cite{isi1925}, originally developed in 1925 to study characteristics of ferro-magnetic material observed in reality, which was applicable to other disciplines like behavior of binary alloys (Bethe \cite{bet1935}). More recent models considering object interaction were built in diverse context. Examples within biology are given by Bramson and Griffeath \cite{bra1980} who examine growth of tumor cells or the model of Kirman \cite{kir1993} explaining herd behavior of ant population in foraging. For information technology state transitions can be used to describe TCP connections or HTTP flows (e.g. Baccelli, McDonald and Lelarge \cite{bac2004} or LeBoudec, McDonald and Mundinger \cite{leb2007}). Haken \cite{hak1983} shows the application to laser light fields produced by excited atoms as well as to chemical and biochemical reactions.
Weidlich \cite{wei1971}) not only examines thermodynamics but also builds the bridge to socio-economics, which is of interest for us. More precisely he interprets the interacting objects as agents of a market with assigned opinions. The interaction and related dynamics are characterized by agents changing their own opinion according to the predominant opinion in the market. Similar models have been developed by many others. The models of F{\"o}llmer \cite{foe1974} as well as Lux \cite{lux1995} might be the most popular.
Although the before mentioned models are well suitable to describe social behavior within a market, they lack a sophisticated link to price dynamics. \\

The microscopic model presented in Pakkanen \cite{pak2010} derives the price dynamics of a single asset by interaction between agents. An agent places a buy or sell order in continuous time which is fulfilled by a so-called market maker holding a sufficient number of shares in order to match supply and demand instantly. Additionally the market maker adjusts the asset price according to the current excess demand. Subsequently, the orders not only impact the price of the asset but also the asset price impacts the agents choice to trade as well as the quantity traded. The mathematical framework on which the model is based is quite interesting as it has a lot of advantages. Compared to other models the author allows for a high degree of individualization related to the agents behavior. The price dynamics are given by a discrete markov chain which is embedded in continuous time using exponentially distributed waiting times leading to a time homogeneous jump-type markov process. In a large market the markov process can then be approximated by a diffusion process. On the other hand, the finite model provides a microscopic foundation for diffusion price processes, heavily used in modern financial mathematics. In order to understand phenomena observable in financial markets, diffusion price processes can be broken down to discrete markov chains, which are often easier to assess analytically. In the primary set-up in Pakkanens model agent interaction is only taking place via feedback through the asset price, which seems unrealistic. Although the model is partially extended in an example by assigning binary opinions to the agents, a general framework for opinion based studies is missing. Moreover the scaling is restrictively chosen to be proportional to the square root of the number of market participants in order to generate stochastic volatility in the diffusion limit.\\ 

In the following we provide a model not only containing a convenient framework for describing socio-economic behavior but also a sophisticated link to price dynamics. Therefore we extend the model elaborated in Pakkanen \cite{pak2010} by providing an additional framework for agent interactions using assigned characteristics. More precisely we assign a state to each agent and measure the endogenous environment by the distribution of all states. We then let the endogenous environment influence each agents tendency to change his state, thus modeling endogenous interaction. We furthermore allow for interaction between the endogenous environment and the price process leading to feedback effects between agents behavior and asset price which are captured as interacting markov chains. In markets with many participants, for the sake of simplicity we want to find a description, which regardless of the number of agents describes the relevant dynamics. To achieve that, the correct scaling of related processes is essential. In contrast to Pakkanen \cite{pak2010} the scaling factor is not fixed in order to be as flexible as possible. See Remark 2.14 in Pakkanen \cite{pak2010} for some more thoughts on the scaling.
Although the following model description is merely an extension of Pakkanen \cite{pak2010}, for better readability we state the extended model rather than referencing to the work done there. Moreover as much as possible the same notation is used.\\

We structure the content as following. In the next section we set the general market framework as a pool of interacting agents trading an asset. To describe their interactive behavior we assign each agent heterogeneous characteristics and express the market as the distribution of those characteristics. The change of characteristics and hence the dynamics of the market is assumed to be stochastic. Then, with defining the agents propensity to trade and specifying the impact on the asset price we link the endogenous market dynamics with the asset price movement. In the third section, we specify conditions under which in a large market the asset price development can be approximated by a diffusion price process. We conclude the section with a theorem summarizing the diffusion approximation and a proposition which appraises the quality of the approximation under specific conditions.
Finally we show the application of the developed model by an example. First we adapt the assumptions made in the model of Lux \cite{lux1995} and derive the endogenous behavior within our model. We then not only compare the resulting price process with the one elaborated in Lux \cite{lux1995} but also provide a large market approximation of the dynamics.
\section{Finite Microscopic Model}

In this section we set the general finite framework of our microscopic model. We extent the model described in Pakkanen \cite{pak2010} by assuming that the asset price depends on another factor, namely the market character which is defined as the distribution of states of the agents participating in the market. With introducing the market character as a key driver of the asset, we can better separate agents behavior and price dynamics. Thus the tractability how the agent, from a rationality and psychological point of view, impacts the price is improved. The rather general term of market character not only includes agents opinion, which makes our model comparable with other opinion based models where the analogous would be the market mood, but also agent type (e.g. noise trader, fundamentalist, guru) or other individual characteristics. We use a markovian framework, that arose out of early models describing phenomena of statistical mechanics and has been the foundation of many models describing interacting objects (see e.g. Kindermann and Snell \cite{kin1980}). Although the markov assumption is rather restrictive, it makes the model memoryless and hence more simple. Additionally the assumption is consistent with the property of diffusion processes heavily used in financial markets.\\
We start with the definition of the endogenous environment by specifying heterogeneous market participants, to which we assign states from a fixed finite set. To reduce complexity we introduce a measure for the distribution of states through the agents leading to the terminology of market character. Then we define the occurrence and  severeness of interaction between the agents, which is modeled as an influence on the decision of a state transition. This results in a dynamical endogenous markovian system in which we measure the related aggregated behavior of the agents with the, now dynamic, market character. Next we specify the individual propensity of the agents to place buy and sell orders and how their actions impact the price. Hereby we explicitly allow for the consideration of external factors in form of random signals. Finally we embed the price process and the market character in continuous time.  We close the section by summarizing the microscopic model and showing the existence of the underlying probability space in a lemma.\\
Let $\mathbb{A}_n = \{1, ..., n\}$ be the \textit{set of agents} participating in the market. To classify common or individual characteristics of the agents we assign to each agent a state, which is used later on to model their endogenous interaction.

\begin{Definition}[State- and configuration space]\ \\
Let $x^a$ be the \textit{state} of agent $a \in \mathbb{A}_n$ which is an element of the fixed finite \textit{state space} $S =\{s_1, ...,  s_m\}, m \in \mathbb{N}$. The vector of all individual states takes values in the compact \textit{configuration space} \mbox{$C := S^{\mathbb{A}_n} = \{x = (x^a)_{a=1}^n, \ x^a \in S \}$}. 
\end{Definition}

During the time $t \in [0, \infty)$ each agent can decide to act (e.g. to change his state). The time of the k-th action is nominated by $T_k \geq 0$, $k \in \mathbb{N}$. The action times are described later in detail, however we use the terminology to describe the development of the states within discrete time in the following definition.

\begin{Definition}[State process]\label{DEF_State_process}\ \\
The state of agent a at time $T_k$ is defined as $x^a_{T_k} \in S$. We capture the development of agent a's state by the process $x_k^a := (x^a_{T_k})_{k \in \mathbb{N}}$ and the development of all agents states by the n-dimensional \textit{state process} $x_k = (x_k^a)_{a=1}^n$. We assume that the vector of initial states is distributed following some n-dimensional distribution function. In particular $x_0 \sim F_{x_0}^n$.
\end{Definition}

In general, the cardinality of the state space $S$ will be much smaller than the one of $\mathbb{A}_n$ (i.e. $m << n$). Moreover, later on we are interested in the development of the market as a whole rather than the development on the level of individual states. Hence it makes sense to coarsen the observable information  for the sake of reduced complexity. 
Rather than the individual states, we consider in the next definition the proportion of all states within the market, representing the overall characteristics of market participants.

\begin{Definition}[Market character]\label{MM}\ \\
For each state we measure the proportion of state $s_i$ among the agents at time $T_k$ by 
\begin{equation}
M_k^i := n^{-d_1} \sum_{a \in \mathbb{A}_n} \mathds{1}_{s_i}(x_k^a), \ k \in \mathbb{N}, \ d_1 \in \mathbb{Q}^+ \geq 1/2.
\end{equation}
The \textit{market character} at time $T_k$ is is defined as the m-dimensional vector valued process of all state proportions, i.e. $M_k = (M_k^i)_{i=1}^m, \ k \in \mathbb{N}$. Additionally, we denote the initial distribution of the market character resulting from Definition \ref{DEF_State_process}, that is the m-dimensional probability distribution of $M_0$ as $F_{M_0}^n$.
\end{Definition}

\begin{Remark}[Scaling of market character]\label{scaling1}
In other models the scaling of the endogenous environment is fixed to be $1/n$ (e.g. Horst and Rothe \cite{hor2008} or Bayraktar, Horst and Sircar et. al \cite{bay2007}) or $1/\sqrt{n}$ (e.g Pakkanen \cite{pak2010}). As concluded in Pakkanen \cite{pak2010}, "Ultimately, the choice of scaling depends on what one wants to model - it seems that $1/\sqrt{n}$ is suited to the study of short-term fluctuations and volatility, whereas $1/n$ is perhaps more appropriate in studies of long-term behavior." We choose a variable scaling factor of $n^{-d_1}, d_1 \in  \mathbb{Q}^+ \geq 1/2$, in order to provide a rather general framework, in which we study can study both. Moreover, note that by construction $n^{d_1-1} M_k$ is a probability measure on the configuration space $C$.
\end{Remark}

\begin{Remark}[Dimension of market character]\label{rem_dim_mm}
While in many opinion-based models only some average-type mood is considered (e.g. Lux \cite{lux1995}, Pakkanen \cite{pak2010}), by construction our market character is m-dimensional. Although adding additional complexity, having a rather general market character provides the necessary flexibility to model agents behavior more specifically. Anyhow, if needed, a reduction of the market character information to relevant properties (e.g. average agent state) is still possible. 
\end{Remark}

We assume that any change of the market is a direct consequence of agents behavior. His behavior is given by so called actions that can either be the change of his state or a trade of the asset. We index each of these actions by $k \in \mathbb{N}$. All information of actions, that have been taken place in the past form the current market history. The k-th action as well as the market history are set more precisely in the next definition.

\begin{Definition}[k-th action, market history]
 \ \\
The \textit{k-th action} is characterized by the tupel $(T_k, A_k, P_k, M_k, B_k), k \in \mathbb{N}$, where $T_k$ is the time when the action occurs,
$A_k \in \mathbb{A}_n$ is the acting agent at time $T_k$ and $B_k \in \{0,1\}$ is an \textit{action indicator} whether the agent trades ($B_k = 1$) or changes his state ($B_k = 0$).  $P_k$ is the \textit{price} per share\footnote{The price is not necessarily be assumed to be logarithmic as in Pakkanen \cite{pak2010}.} and $M_k$ the above mentioned character of the market. All information is captured in the \textit{market history}, which is given by $\mathcal{G}_k := \sigma (T_i, A_i, P_i,M_i,B_i , i \leq k)$. 
\end{Definition}

\begin{Assumption}[k-th action]\label{assumption_k-th_action}\ \\
We require that only one agent is acting at a specific point in time as well as that the acting agent either trades or changes his state. Although the first part of the assumption seems rather strong, it is however reasonable as actions are performed in continuous time and are very unlikely to happen at the same point in time. The dichotomous behavior of any agent is assumed mainly for the reason of simplicity, as it leaves the market character and the price process rather separable.
\end{Assumption}

Next, we specify the tendency of each agent to act before characterizing the action and related impact.
To determine the likelihood of agent a to be the one who acts at time $T_k$, we assign to each agent intensity (or rate-) functions and then weight the agents. \\
The agent specific tendency to act (i.e. to trade or to switch his state) is assumed to be dependent on the price as well as the character of the market. For the propensity to trade we use the trading intensity function defined in Pakkanen \cite{pak2010}, but allow additionally for dependence on the character of the market.  
The state transition rate function is defined analogously, depending on both, market character and price. We assume that there is always an agents who wants to act and thus require the trading intensity function as well as the transition rate function to be positive.

\begin{Definition}[Trading intensity, state transition rate, action rate]\label{def_intensity} \ \\
Let $ \lambda_a \colon \mathbb{R} \times \mathbb{R}^{m} \to \mathbb{R}_+, a \in \mathbb{A}_n $  be the continuous and bounded \textit{trading intensity function} of agent $a$.
Moreover the \textit{aggregate trading intensity} is defined as the sum of trading intensities of over all agents via $\lambda_{\mathbb{A}_n} := \sum^n_{a=1} \lambda_a$. \\ 
Similarly let $ \mu_a \colon \mathbb{R} \times \mathbb{R}^{m} \to \mathbb{R}_+, a \in \mathbb{A}_n $ be the \textit{state transition rate function}, which is assumed to be continuous and bounded and denote the \textit{aggregate state transition rate} by  $\mu_{\mathbb{A}_n} := \sum^n_{a=1} \mu_a$. We summarize the intensity of all actions with the \textit{aggregated action rate} $\nu_{\mathbb{A}_n}(x,v) := \lambda_{\mathbb{A}_n}(x,v) + \mu_{\mathbb{A}_n}(x,v)$. 
\end{Definition}

In the next definition we specify the acting probabilities of individual agents. Heuristically we weight the respective intensity or rate function.	

\begin{Definition}(Acting probabilities)\ \\
The probability, that agent a trades at $T_k$ is defined as 
\begin{equation}\label{prob_trade}
\p (A_k = a, B_k = 1|\mathcal{G}_{k-1}) = \frac{\lambda_a (P_{k-1},M_{k-1}) }{\nu_{\mathbb{A}_n} (P_{k-1},M_{k-1})}.
\end{equation}

Similarly, we define the probability, that agent a changes his state by

\begin{equation}\label{prob_trade2}
\p (A_k = a, B_k = 0|\mathcal{G}_{k-1}) =\frac{\mu_a (P_{k-1},M_{k-1}) }{\nu_{\mathbb{A}_n} (P_{k-1},M_{k-1})}.
\end{equation}

Moreover the probability that the k-th action is a state transition is set as

\begin{equation}
\p (B_k = 0|\mathcal{G}_{k-1}) = \sum_{a=1}^n \frac{\mu_a (P_{k-1},M_{k-1}) }{\nu_{\mathbb{A}_n} (P_{k-1},M_{k-1})} = \frac{\mu_{\mathbb{A}_n} (P_{k-1},M_{k-1}) }{\nu_{\mathbb{A}_n} (P_{k-1},M_{k-1})},
\end{equation}

and analogously the probability that the k-th action is a trade is given by

\begin{equation}
\p (B_k = 1|\mathcal{G}_{k-1}) = \frac{\lambda_{\mathbb{A}_n} (P_{k-1},M_{k-1}) }{\nu_{\mathbb{A}_n} (P_{k-1},M_{k-1})} = 1- \p (B_k = 0|\mathcal{G}_{k-1}).
\end{equation}

\end{Definition}
The next step is to characterize the state transition laws and consequential derive the dynamics of the market character. Although the probability is not further determined here, we introduce an extra notation to clarify that we explicitly allow for dependence of individual state transition probabilities on the market character and price.   

\begin{Definition}[State transition probability]\label{def_sp}\ \\
We use the following notation for the individual state transition probability, i.e. the probability that agent a changes from $s_i$ to $s_j$, given that he is the one that wants to change his state.
\begin{equation}
\Pi^{i,j}_{n,a}(P_{k-1},M_{k-1}) := \p (x_k^a = s_j|x_{k-1}^a = s_i, B_k = 0, A_k = a)
\end{equation}
\end{Definition}

We capture all state transition probabilities per agent in a transition matrix, i.e. we define
\begin{equation}
\Pi_{n,a} := (\Pi^{i,j}_{n,a})^m_{i, j = 1}
\end{equation}

While Definition \ref{def_sp} is quantifying a single movement from state $s_i$ to $s_j$, we are rather interested in the aggregated dynamics. The aggregated behavior of all agents, i.e. are agents rather joining or leaving a state is used to describe attractiveness of a state.

\begin{Definition}[Aggregated state transition]\ \\
Let $\Pi^{i-}_{n,a}$ describe the aggregated propensity to leave state $s_i$, i.e. the probability of a state transition from $s_i$ to any other state at time $T_k$ More precisely,
\begin{equation}\label{aggregate_minus}
 \Pi^{i-}_{n,a} (P_{k-1},M_{k-1}) := n^{d_1-1} M_{k-1}^i \sum^m_{\substack{j = 1, \\ j \neq i}} \Pi_{n,a}^{i,j}(P_{k-1},M_{k-1}),
\end{equation}

where the pre-factor $n^{d_1-1} M_{k-1}^i$ is the probability that the acting agent $A_k$ had state $s_i$, which is well defined as $n^{d_1-1} M_{k-1}$ is a probability measure on $C$.
Analogously define the aggregated propensity to switch to state $s_i$ by 

\begin{equation}\label{aggregate_plus}
 \Pi^{i+}_{n,a} (P_{k-1},M_{k-1}) := n^{d_1-1} \sum^m_{\substack{j = 1, \\ j \neq i}} M_{k-1}^j \Pi_{n,a}^{j,i}(P_{k-1},M_{k-1}). 
\end{equation}

\end{Definition}

We are now able to derive the dynamics of the market character in the following lemma. As we assume that only one agent can act on each action time $T_k$ the proportion of a state at time $T_k$ can either increase or decrease by $n^{-d_1}$ or stay unchanged. 

\begin{Lemma}[Market character dynamics]\label{MMDynamics}\ \\
The probability that an agent of state $s_i$ switches to a different state $s_j$ and therefore that the proportion of state $s_i$ decreases by $n^{-d_1}$ is given by 
\begin{equation}\label{dynamics_mood2}
\begin{split}
\p (M_k^i - M_{k-1}^i = -n^{-d_1} | \mathcal{G}_{k-1}) 
&= \frac{\sum^n_{a = 1} \mu_{a} (P_{k-1},M_{k-1}) \Pi^{i-}_{n,a} (P_{k-1},M_{k-1})}{\nu_{\mathbb{A}_n} (P_{k-1},M_{k-1})} 
\end{split}
\end{equation}

and similarly the probability that the occupancy measure increases by $n^{-d_1}$ is given by
\begin{equation}\label{dynamics_mood3}
\begin{split}
\p (M_k^i - M_{k-1}^i = n^{-d_1} | \mathcal{G}_{k-1}) 
&= \frac{\sum^n_{a = 1} \mu_{a} (P_{k-1},M_{k-1}) \Pi^{i+}_{n,a} (P_{k-1},M_{k-1})}{\nu_{\mathbb{A}_n} (P_{k-1},M_{k-1})}.
\end{split}
\end{equation}
\begin{proof}
Using the fact that $n^{d_1-1} M_{k-1}^i$ is a probability measure on $C$ as well as the representation defined in Equation (\ref{prob_trade2}) and (\ref{def_sp}) the lemma follows from basic calculations.
\end{proof}
\end{Lemma}

After we set the general framework for the endogenous environment we now can describe the impact on the price and feedback effects. Therefore we specify the discrete price chain $(P_k)_{k=0}^\infty$, the dynamics of which we want to be dependent on endogenous as  well as exogenous factors. Although the main part of the model consists of endogenous factors (i.e. market character and feedback effects of the price) in order to be flexible as well as to be consistent with Pakkanen \cite{pak2010}, we allow the price dynamics to be dependent on a random signal describing exogenous impacts. 
The volume of shares agent $A_k$ would trade at $T_k$ is quantified by the \textit{excess demand function} defined next.

\begin{Definition}[Excess demand functions, signals]\ \\
The \textit{excess demand function} $e^n_a \colon \mathbb{R}^{m+2} \to \mathbb{R}, \ a \in \mathbb{A}_n$ is a measurable function depending on $P_{k-1}$, $M_{k-1}$ and the random variable $\xi_k$, which is assumed to be independent of $\mathcal{G}_{k-1}$ and $A_k$. Additionally, we assume that the \textit{signals} $(\xi_k)^\infty_{k=1}$ are i.i.d. with cdf $F_\xi$.
\end{Definition}

The price at time $T_k$ will be set by a market maker, which is assumed to handle all trades, and is defined by a pricing rule depending on the excess demand of the acting agent and the old price.

\begin{Definition}[Pricing rule]\ \\
Consider the borel measurable \textit{pricing rule}\footnote{As we could factor out the market character to the excess demand function, where it is already considered, we refrain from including it in the pricing rule.} function $r_n \colon \mathbb{R}^2 \to \mathbb{R}$ setting the price $P_k$ via
\begin{equation}
P_k = r_n(e^n_{A_k}(P_{k-1}, M_{k-1}, \xi_k), P_{k-1}).
\end{equation} 
\end{Definition}

By construction $(P_k)_{k=0}^\infty$ and $(M_k)_{k=0}^\infty$ now are two interacting Markov chains. In order to embed them homogeneously in continuous time and thus describing the price as well as the character by a time homogeneous Markov process, we further characterize the points in times in which the agents decide to act. 

\begin{Definition}[Intra-action times]\label{Def_intra-action}\ \\
The \textit{intra-action times} $(\tau_k)^\infty_{k=1}$ are defined as $\tau_k := T_k - T_{k-1}, k \geq 1.$\\
Since we want the intra-action times to be memory-less for the sake of simplicity, i.e.
\begin{equation}
\p (\tau_k > t +h | \tau_k > h, \mathcal{G}_{k-1}) = \p(\tau_k > t | \mathcal{G}_{k-1}), \ t,h \geq 0.
\end{equation}
the intra-action times are assumed to be exponentially distributed. Heuristically we assume that the rate of the exponential distribution is given by the aggregated action rate, i.e.
\begin{equation}\label{exponential_waiting_times}
\p (\tau_k \in [0,t]|\mathcal{G}_{k-1}) = 1 - e^{-\nu_{\mathbb{A}_n}(P_{k-1},M_{k-1}) t}, \ t \geq 0,
\end{equation}
More precisely, to ensure a sufficient level of independence between the source of randomization and the price as well as market character we need to assume that the Intra-action times $(\tau_k)^\infty_{k=1}$ are given by 
\begin{equation}
\tau_k := \frac{\gamma_k}{\nu_{\mathbb{A}_n}(P_{k-1},M_{k-1})}, \ k \in \mathbb{N},
\end{equation} 
where $(\gamma_k)^\infty_{k=1}$ are i.i.d. random variables independent of $(P_k,M_k)^\infty_{k=0}$ with $\gamma_1 \sim Exp(1)$.
\end{Definition}

\begin{Definition}[Price process, market character index]\label{def_priceprocess}\ \\
After setting an initial price $P_0 \sim F_{P_0}$ and fixing $T_0 = 0$ we can define the \textit{price process} as
\begin{equation}
X^n_t := \sum^\infty_{k=0} P_k \mathds{1}_{[T_k,T_{k+1})}(t), \ t \geq 0.
\end{equation} 

Analogously we introduce the \textit{market character index} via
\begin{equation}
V^n_t := \sum^\infty_{k=0} M_k \mathds{1}_{[T_k,T_{k+1})}(t), \ t \geq 0.
\end{equation} 
\end{Definition}

Note that by construction $X^n_t$ and $V^n_t$ are c\'{a}dl\'{a}g and that $F^n_{M_0}$, in contrast to $F_{P_0}$, is depending on $n$. 

The next lemma now summarizes the construction of the finite microscopic model. It states the existence of a probability space carrying the price process as well as the market character index as time homogeneous Markov processes. Furthermore it gives the rate kernel as the product of action rate and transition kernel. The basis of the lemma builds the synthesis theorem (e.g. Theorem 12.18 of Kallenberg \cite{kal2002}), which embeds a discrete Markov chain into continuous time using exponentially distributed waiting times.  

\begin{Lemma}[Existence]\label{lemma_existence}\ \\
If the preceding Assumptions hold, then there exists a probability space $(\Omega, \mathcal{F}, \mathbb{P})$ which carries the model such that $(X^n_t,V^n_t)_{t \in [0, \infty)}$ is a time homogeneous pure jump Markov process with rate kernel 

\begin{equation}\label{rate_kernel}
K_n(x,v,dy,dw) := \nu_{\mathbb{A}_n}(x,v)  k_n(x,v,dy,dw),
\end{equation}

where the transition kernel $k_n(x,v,dy,dw)$ is a regular version of the conditional distribution $\mathbb{P}(P_1 - P_0 \in dy, M_1 - M_0 \in dw | P_0 = x, M_0 = v)$.
\begin{proof}
By the construction of the Markov chain $(P_k,M_k)_{k=0}^\infty$ and the assumption made in Definition \ref{Def_intra-action} the synthesis theorem (e.g. Theorem 12.18 of Kallenberg \cite{kal2002}) states that $(X^n_t,V^n_t)_{t \in [0, \infty)}$ is a pure jump-type Markov process and also gives the rate kernel. Time homogeneity is given by the recursive definition of $(P_k,M_k)_{k=0}^\infty$ (see e.g. Proposition 8.6 of Kallenberg \cite{kal2002}) as the pricing rule $r_n$ as well as the transition matrix $\Pi_{n,a}$ are independent of the time.
\end{proof}
\end{Lemma}

\section{Diffusion approximation}\label{sec_diff_approx}

In this chapter we have a closer look on markets with a large number of participants. For the sake of simplicity we want to describe the relevant dynamics regardless of the number of market participants. We employ (It{\^o}-) diffusion processes, which are commonly used in financial mathematics to model financial instruments. As such we not only have a powerful tool to approximate the dynamics of a large market within our model, but also use our microscopic model to explain foundations of various extraordinary characteristics observed in the price processes in real market data. One example - phase transitions and oscillations - is discussed in the second part of this paper.
In the following, we give the general conditions under which the indices of our microscopic model can be approximated by a (It{\^o}-) diffusion process $(X_t,V_t)_{t \in [0,\infty)}$, which is described as the solution of the $m+1$ dimensional stochastic differential equation (SDE)
\begin{equation}
d(X_t,V_t) = \hat{b}(X_t,V_t)dt + \hat{\sigma}(X_t,V_t)dB_t, 
\end{equation} 
where $(X_t,V_t) \in \mathbb{R} \times \mathbb{R}^m$ and $B_t$ is a $m+1$ dimensional brownian motion. In the main theorem of this chapter we not only proof the existence of a limit of $(X^n_t,V^n_t)_{t \in [0,\infty)}$ as $n \to \infty$, but also identify the coefficients $\hat{b}$ and $\hat{\sigma}$ as the limit of the first, respectively second, moment of the finite process $(X^n_t,V^n_t)_{t \in [0,\infty)}$.

To ensure the convergence of the time-homogeneous pure jump process $(X^n_t,V^n_t)_{t \in [0,\infty)}$ to a (continuous) diffusion as $n \to \infty$, we need to make several additional assumptions. Firstly we need to ensure the existence and convergence of the first and second moments of the price process $X_t^n$ as well as of the market character index $V_t^n$. 
Secondly, the moments should meet some regularity conditions in the limit. Although convergence to a continuous diffusion can be achieved under different regularity conditions (see e.g. Mao \cite{mao2007} or Xua et al. \cite{xua2008} for more details) we choose the locally lipschitz and linear growth condition to be consistent with Pakkanen \cite{pak2010}.
Finally, we want the action rate function $\nu_{\mathbb{A}_n}$ to converge "nicely" in order to not have jumps in the limit of $(X^n_t,V^n_t)$. \\
To characterize the first moment of $X_t^n$ we introduce a function which quantifies the aggregated expectation in terms of demand and supply in our market in the next definition. The so-called expected aggregate excess demand $z_n$ is the individual expected excess demand of the agents aggregated by weighting their respective trading intensity function and depends on the market character as well as on the price.
%

\begin{Definition}[Excess demands and pricing rule]\label{assumption_excess_demands}\ \\
Let the \textit{expected aggregate excess demand} at $(x,v) \in \mathbb{R} \times \mathbb{R}^{m}$ be defined as
\begin{equation}\label{z_n}
z_n(x,v) := n^{-d_2}  \sum^n_{a=1} \lambda_a(x,v) \E [e^n_a(x,v,\xi_1)], \ d_2 \in \mathbb{Q}^+ \geq 1/2.
\end{equation}

Furthermore we assume the pricing rule $r_n$ is given by
\begin{equation}\label{pricing_rule}
r_n(q,x) = x + \alpha n^{-d_2} q + u_n(q,x), \ q, x \in \mathbb{R}, \ d_2 \in \mathbb{Q}^+ \geq 1/2,
\end{equation}
where $\alpha > 0$ and $u_n, n \in \mathbb{N}$ is a borel measurable function such that $\forall \delta >0 \ \exists C_n^\delta, \ n \in \mathbb{N}$ such that $C_n^\delta = o(n^{-1})$ and $\sup_{|x| \leq \delta} |u_n(q,x)| \leq C_n^\delta |q|, \ \forall q \in \mathbb{R}, \ n \in \mathbb{N}$.
\end{Definition}

As visible in Equation (\ref{pricing_rule}) we assume the pricing rule to be nearly affine, that is affine apart from a function $u_n$ which is bounded by a constant $C_n^\delta$ converging to zero when the number of agents tents to infinity. So we ensure that a possible "nice" behavior of excess demands is sufficiently carried to the increments of $X_t^n$. Thereby we allow for a flexible scaling factor $n^{-d_2}$ following the same rational as stated in Remark \ref{scaling1}.\\
 
Analogously to Equation (\ref{z_n}) we define the expected aggregated transition related to state $s_i$ by summing the individual agent transitions weighted by the respective transition rate. The resulting expected aggregate state transition then describes the first moments of $V_t^n$.

\begin{Definition}[Expected aggregated state transition]\label{assumption_state_movement}\ \\
We define the \textit{expected aggregate state transition of state $s_i$} as
\begin{equation}\label{b_n^i}
b_n^i(x,v) := n^{-d_1} \sum^n_{a=1} \mu_a(x,v)\left( \Pi_{n,a}^{i+}(x,v) - \Pi_{n,a}^{i-}(x,v) \right).
\end{equation} 
In summary, we write the \textit{expected aggregate state transition} as

\begin{equation}\label{b_n}
b_n(x,v) := 
\begin{bmatrix} 
b_n^1\\ 
b_n^2\\
\vdots\\
b_n^m\\
 \end{bmatrix} (x,v).
\end{equation}

\end{Definition}

In the following Assumption we quantify the second moments of the price process. In particular we define the \textit{expected trading volume} by aggregating the second moments of the heterogeneous excess demand functions weighted by the respective intensity function.

\begin{Definition}[Trading volume]\label{assumption_trading_volume}\ \\
The \textit{expected trading volume} at price and character level $(x,v) \in \mathbb{R} \times \mathbb{R}^m$ is defined by
\begin{equation}
\sigma_n(x,v) := \left( n^{-2d_2} \sum^n_{a=1} \lambda_a(x,v) \E [e^n_a(x,v,\xi_1)^2]\right)^{1/2}
\end{equation}

\end{Definition}

Analogous to Definition \ref{assumption_trading_volume} we describe the second moment of the market character by aggregating the second moments of the individual state transition weighted by the respective rate function. The resulting so-called \textit{transition volume} is described by the variance within states $s_i$ on the one hand and covariances between states $s_i$ and $s_j$ on the other.

\begin{Definition}[Transition volume]\label{var_states}\ \\
We denote the \textit{transition volume between $s_i$ and $s_j$} with
\begin{equation}
c_n^{i,j}(x,v) := \left( -n^{-(d_1+1)} \sum^n_{a=1} \mu_a(x,v) (v_i \Pi^{i,j}_{n,a}(x,v) + v_j \Pi^{j,i}_{n,a}(x,v)) \right)^{1/2}, \ (x,v) \in \mathbb{R} \times \mathbb{R}^m
\end{equation}
and the \textit{transition volume within $s_i$}
\begin{equation}\label{c_n^i}
c_n^{i}(x,v) := \left( n^{-2d_1} \sum^n_{a=1} \mu_a(x,v) (\Pi^{i +}_{n,a}(x,v) + \Pi^{i -}_{n,a}(x,v)) \right)^{1/2}, \ (x,v) \in \mathbb{R} \times \mathbb{R}^m.
\end{equation}
In short, we write
\begin{equation}\label{c_n}
c_n(x,v) := 
\begin{bmatrix} 
c_n^1 & c_n^{1,2} & ... & c_n^{1,m}\\ 
c_n^{2,1}   & c_n^2 & \ & \vdots \\
\vdots & \ & \ddots & \ \\
c_n^{m,1} & ... & \ & c_n^m\\
 \end{bmatrix} (x,v)
\end{equation}

and call the function $c_n$ \textit{transition volume}. 
\end{Definition}

In order to achieve convergence of $(X^n_t,V^n_t)_{t \in [0,\infty)}$ to a continuous diffusion neither jump size nor the intensity should explode. While the jump size of $V_t^n$ is bounded by the construction of the difference $M_k - M_{k-1}$, we need restrictions on the excess demands in order to bound the jump size of $X_t^n$. We bound the trading intensity as well as the state transition rate by restricting the action rate.

\begin{Assumption}[No explosions]\label{assumption_explosion}\ \\
For every $\delta > 0$,
\begin{enumerate}
\item  $\limsup_{n \to \infty} \sup_{|(x,v)| \leq \delta}|\frac{\nu_{\mathbb{A}_n}(x,v)}{n}| < \infty$ and
\item  $\left\{ e^n_a(x,v,\xi_1)^2 \colon |(x,v)| \leq \delta, \ a \in \mathbb{A}_n, \ n \in \mathbb{N} \right\}$ is uniformly integrable
\end{enumerate}
\end{Assumption}

Now, we are in the position to apply Theorem IX. 4.21 of Jacod and Shiryaev \cite{jac2003}, which gives the convergence of the process $(X^n_t,V^n_t)_{t \in [0,\infty)}$ to a diffusion process when the number of market participants tends to infinity.
The drift coefficient is determined by functions $z$ and $b$ defined in Definition \ref{assumption_excess_demands} and \ref{assumption_state_movement}, while the diffusion coefficient is given by functions $\sigma$ and $c$ described in Definition \ref{assumption_trading_volume} and \ref{var_states}. We summarize the diffusion approximation in the following theorem.

\begin{Theorem}[Diffusion approximation]\label{Diff_approx}\ \\
Assume that for the functions $z_n, b_n, \sigma_n$ and $c_n$ there exist continuous functions $z, b, \sigma$ and $c$ that are locally lipschitz and of linear growth such that $z_n \to z, b_n \to b, \sigma_n \to \sigma$ and $c_n \to c$ uniformly on compact sets (u.o.c.) for $n \to \infty$. If additionally Assumption \ref{assumption_explosion} holds and $F_{M_0}^n \xrightarrow{n \to \infty} F_{M_0}$, then

\begin{equation}
(X^n_t, V^n_t)_{t \in [0, \infty)} \xrightarrow{\mathcal{L}} (X_t,V_t)_{t \in [0, \infty)} \ in \ D_{\mathbb{R}^{m+1}}[0, \infty),
\end{equation}

where $(X_t,V_t)_{t \in [0, \infty)}$ is the unique strong solution of the SDEs

\begin{equation}\label{SDE1}
\begin{cases}
dX_t = \alpha z (X_t, V_t)dt + \alpha \sigma(X_t,V_t) dB_t, & X_0 = \zeta\\
dV_t = b(X_t,V_t)dt + c(X_t,V_t) d\textbf{B}_t, & V_0 = \theta ,
\end{cases}
\end{equation}

where $(B_t)_{t \in [0, \infty)}$ is a one dimensional standard Brownian motion, $\zeta \sim F_{P_0}$ independent of $B_t$, and $(\textbf{B}_t)_{t \in [0, \infty)}$ is a $m$-dimensional Brownian motion, which is independent of $\theta \sim F_{M_0}$.
\begin{proof}
See Appendix \ref{proof_diff_approx1}.
\end{proof}
\end{Theorem}

In case the large market limit of the price- and market character index given by Theorem \ref{Diff_approx} is deterministic (i.e. $\sigma_n = 0$ and $c_n = 0$), the rate by which the pure-jump type process $(X^n_t, V^n_t)_{t \in [0, \infty)}$ converges to the limit process $(X_t,V_t)_{t \in [0, \infty)}$ when $n \to \infty$ can be assessed. The following Proposition gives particularly the convergence rate as being the speed by which $\sigma_n$, respectively $c_n$, tend to zero. 

\begin{Proposition}[Rate of convergence]\label{prop_rate_of_conv}
Assume that $F_{M_0}^n = F_{M_0} $. Let $(a_n)_{n \geq 0}$ be a positive sequence with $a_n \to \infty$ such that
\begin{enumerate}
\item $a_n^2 \sigma_n^2 \xrightarrow[n \to \infty]{\text{u.o.c.}} \hat{\sigma}^2$ and $a_n^2 c_n^2 \xrightarrow[n \to \infty]{\text{u.o.c.}} \hat{c}^2$ for some continuous functions $\hat{\sigma}, \hat{c}$
\item $\sqrt{n} a_n = O(n^{d_1}+n^{d_2}+1/C_n^{\delta})$
\item $a_n(z_n-z) \xrightarrow[n \to \infty]{\text{u.o.c.}} 0$ and $a_n(b_n-b) \xrightarrow[n \to \infty]{\text{u.o.c.}} 0$.
\end{enumerate}

Then 
\begin{equation}
\sup_{s \leq t} |(X^n_s,V^n_s)-(X_s,V_s)| \leq a_n^{-1} \sup_{s \leq t}|(Y_s,Z_s)|, \ \forall t \geq 0,
\end{equation}
where $(Y_t,Z_t)_{t \in [0, \infty)}$ is the solution of the SDE

\begin{equation}
\begin{cases}
dY_t = \hat{\sigma}(Y_t, Z_t)dB_t, & Y_0 = 0\\
dZ_t = \hat{c}(Y_t,Z_t) d\textbf{B}_t, & Z_0 = 0.
\end{cases}
\end{equation}

\begin{proof}
See Appendix \ref{proof_prop_rate_of_conv}.
\end{proof}
\end{Proposition}

\section{Example: Phase transitions and oscillations}

Inspired by an observation in entomology, in particular related to ant populations and their contagious behavior towards food collection as discussed in Kirman \cite{kir1993}, Lux  \cite{lux1995} applied a herding mechanism to a fixed number of noise traders to describe an asset market. Thereby the individual noise trader is either optimistic or pessimistic and the rule for opinion change depends on the opinion of the majority as well as price trends. In order to have price changes be determined by a market maker through supply and demand matching, Lux introduces a second type of traders, so called fundamentalists that sell (buy) when the price is above (below) a fundamental value. Lux then used the master equation approach, originated from elementary particle systems in physics (see e.g. Haken \cite{hak1983}), together with methods discussed in Weidlich and Haag \cite{wei1983} to derive the properties of his market and to show its capability of generating bubbles and periodic oscillations. 

In this chapter we embed the model described in Lux \cite{lux1995} within our framework. We model the endogenous behavior of the agents according to the assumptions made in the model of Lux and assess if we achieve similar results, i.e herd behavior leads to temporary equilibria of the proportion of optimists and pessimists. We then link the endogenous dynamics with a price process and show its convergence to a diffusion when the market is large.
This section is structured as following. In the first sub-section we establish the endogenous environment using only noise traders, which builds the base of the market character index. We then derive its properties and compare the results to Lux \cite{lux1995}. In a second step we discuss the diffusive limit of the market character index and assess the rate of convergence. In the second sub-section we introduce fundamentalists as an additional group of traders and link the endogenous environment with the price process and vice versa. Hereby, in contrast to Lux \cite{lux1995} we allow additionally for random signals influencing the agents trade behavior. We again assess the properties, derive the diffusion approximation and make a concluding comparison with Lux \cite{lux1995}.

\subsection{Pure endogenous dynamics}
\subsubsection{Finite Model}
First we specify all model components that directly affect the market character. We start with a finite set of an even\footnote{We chose the number of agents to be even in order to be consistent with Lux \cite{lux1995}.} number of agents $\mathbb{A}_n = \{1, ..., n\}, n=2 \mathbb{N}$. Moreover we consider a state space $S = \{-1,1\}$ where $s_1 = -1$ represents a pessimistic and $s_2 = 1$ an optimistic opinion. Since no further specifications of the distribution of initial states is given in Lux \cite{lux1995}, we assign each agent $a \in \mathbb{A}_n$ an initial state $x^a_0 \in S :=\{-1,1\}$ such that the vector of initial states $x_0$ has some probability distribution $F_{x_0}^n$.
In line with Lux \cite{lux1995} we choose the scaling of the market character  to be $\frac{1}{n}$, i.e $d_1=1$. Following these assumptions, the character of the market at time $T_k $ is given by $M_k = (M_k^1,M_k^2)$ with 
\begin{equation}\label{M_k^1_Lux}
M_k^1 = \frac{1}{n} \sum_{a \in \mathbb{A}_n} \mathds{1}_{\{-1\}}(x_k^a), \ k \geq 0,
\end{equation}

\begin{equation}\label{M_k^2_Lux}
M_k^2 = \frac{1}{n} \sum_{a \in \mathbb{A}_n} \mathds{1}_{\{1\}}(x_k^a), \ k \geq 0.
\end{equation}

Since $n$ is constant and the market is dichotomous, the two-dimensional market character is fully described by the one-dimensional \textit{average opinion} defined as
\begin{equation}\label{Def_overline_M}
\overline{M}_k := \frac{1}{n} \sum_{a \in \mathbb{A}_n} x_k^a, \ k \geq 0
\end{equation}
by $M_k = \left(\frac{1-\overline{M}_k}{2},\frac{1+\overline{M}_k}{2}\right)$.

Considering the average opinion not only reduces the dimension and thus simplifies the endogenous dynamics, but also is consistent with the \textit{average opinion} examined in Lux \cite{lux1995}. We denote the initial distribution of $\overline{M}_0$ which results from $F_{x_0}^n$  by Equation (\ref{Def_overline_M}) by $F_{\overline{M}_0}^n$.

Next we define the state transition probability, which describes the likelihood an optimistic, respective pessimistic, agent is to change his opinion. Note that in line with the homogeneity assumption made in Lux \cite{lux1995} the transition probabilities are common to all agents.

\begin{Definition}[Transition probabilities]\ \\
The transition probability to switch the state from $-1$ to $1$ is defined as
\begin{equation}\label{state1_example_lux}
\Pi^{1,2}(\overline{M}_{k-1}) = \beta e^{\gamma \overline{M}_{k-1}}
\end{equation}
and analogous the probability to switch the state from $1$ to $-1$ is defined as
\begin{equation}\label{state2_example_lux}
\Pi^{2,1}(\overline{M}_{k-1}) = \beta e^{-\gamma \overline{M}_{k-1}},
\end{equation}
where $\beta, \gamma >0$ and $\beta < e^{-\gamma}$ and hence the transition matrix, which is assumed to be common to all agents, is given by
\begin{equation}
\Pi_n(\overline{M}_{k-1}) = \begin{pmatrix}  \Pi^{1,2} & 1-\Pi^{1,2}   \\ 1-\Pi^{2,1}   & \Pi^{2,1} \\ \end{pmatrix} (\overline{M}_{k-1}).
\end{equation}
\end{Definition}

\begin{Remark}
The explicit form of transition probabilities presented in Equations (\ref{state1_example_lux}) and (\ref{state2_example_lux})  were chosen by Lux to reflect the following socioeconomic characteristics. Firstly, the transition probability needs to reflect the idea of herding, i.e. the tendency of an agent to change his opinion to be optimistic (pessimistic) is larger when the majority of the traders already has an optimistic (pessimistic) opinion. Moreover the relative change in probability should change linear with the majority's opinion and be symmetric for optimism and pessimism, viz. 
\begin{equation}
\frac{\partial\Pi^{1,2}(\bar{v})}{\Pi^{1,2}(\bar{v})} = C d\bar{v} = -\frac{\partial\Pi^{2,1}(\bar{v})}{\Pi^{2,1}(\bar{v})}, \bar{v} \in [-1,1],
\end{equation}
for some constant $C \neq 0$. Finally, by definition, the probability needs to be between zero and one. The functional form of Equations (\ref{state1_example_lux}) and (\ref{state2_example_lux}) not only meets the requirements above but also give a good control of the infection by the parameters $\beta$ and $\gamma$. While $\gamma$ regulates the intensity of the infection and thus herd behavior, $\beta$ controls the speed of contagion and hence contributes to the time scale.
\end{Remark}

Since the agents are assumed to behave homogeneously we assume a common state transition rate $\mu_a$. Moreover as we are rather interested in the agent interaction itself and less on the time scale and could factor the transition rate into $\beta$ anyway, we set 

\begin{equation}\label{mu_a_lux_endo}
\mu_a = 1, \forall a \in \mathbb{A}_n.
\end{equation}

Since the average opinion $\overline{M}_k$ can from time $T_k$ to $T_{k+1}$ either change by $\pm \frac{2}{n}$ or stay unchanged, it has its values on the $n+1$ valued lattice $\mathbb{L}$ from -1 to 1, viz.
\begin{equation}
\overline{M}_k \in \mathbb{L}, \ \forall k \geq 0, \ \text{with} \ \mathbb{L}:= \left\{ -1, -\frac{n-2}{n}, -\frac{n-4}{n}, \dots, \frac{n-4}{n}, \frac{n-2}{n}, 1  \right\}.
\end{equation}

In summary, $(\overline{M}_k)_{k=0}^{\infty}$ is a Markov chain on $\mathbb{L}$ with state dependent transition probabilities, which are by Lemma \ref{MMDynamics} given as

\begin{equation}\label{overline_M_prob1}
\begin{split}
\p \left( \overline{M}_{k} - \overline{M}_{k-1} = \frac{2}{n} | \mathcal{G}_{k-1} \right) &= \p \left( M^2_k - M^2_{k-1} =  \frac{1}{n} | \mathcal{G}_{k-1} \right)\\
&= M^1_{k-1} \Pi^{1,2}(\overline{M}_{k-1})\\
&= \frac{1-\overline{M}_{k-1}}{2} \Pi^{1,2}(\overline{M}_{k-1})
\end{split}
\end{equation}
and 
\begin{equation}\label{overline_M_prob2}
\begin{split}
\p \left( \overline{M}_{k} - \overline{M}_{k-1} = -\frac{2}{n} | \mathcal{G}_{k-1} \right) &= \p \left( M^2_k - M^2_{k-1} =  -\frac{1}{n} | \mathcal{G}_{k-1}\right)\\
&= M^2_{k-1} \Pi^{2,1}(\overline{M}_{k-1})\\
&= \frac{1+\overline{M}_{k-1}}{2} \Pi^{2,1}(\overline{M}_{k-1})
\end{split}
\end{equation}

Following Definition \ref{def_priceprocess} we embed the markov chain $(\overline{M}_k)_{k \geq 0}$ in continuous time using the \textit{average opinion index}, defined as

\begin{equation}\label{def_average_mood_process}
\overline{V}^n_t := \sum_{k=0}^{\infty} \overline{M}_k \mathds{1}_{[T_{k},T_{k+1})}(t), \ t \geq 0. 
\end{equation}

Note that by Equation (\ref{exponential_waiting_times}) and (\ref{mu_a_lux_endo}) we have for the intra-action times: $\tau_k \sim Exp(n)$.
In order to study the stationary behavior of $\overline{V}^n_t$, we calculate the stationary distribution of the underlying Markov chain $\overline{M}_k$.

\begin{Proposition}[Stationary distribution]\label{prop_stat_dist}\ \\
The stationary distribution of $\overline{M}_{k}$ and $\overline{V}_{t}^n$ resulting from Equations (\ref{overline_M_prob1}), (\ref{overline_M_prob2}) and (\ref{def_average_mood_process}) is given by 
\begin{equation}\label{stationary_dist}
\p_{st}(\bar{v}) = \p_{st}(0) \frac{\left(\frac{n}{2} !\right)^2}{n!} \binom{n}{\frac{n(1+\bar{v})}{2}} \exp\left(\frac{\gamma n \bar{v}^2}{2}\right), \bar{v} \in \mathbb{L},
\end{equation}
where $\p_{st}(0)$ is determined by the normalization condition
\begin{equation}
\sum_{\bar{v} \in \mathbb{L}} \p_{st}(\bar{v}) = 1.
\end{equation}
If
\begin{equation}
\gamma \underset{(<)}{>}\frac{n}{2} \ln \left(\frac{n+2}{n}\right),
\end{equation}
then $\p_{st}$ has a local minimum (maximum) at 0.

\begin{proof}
See Appendix \ref{proof_prop_stat_dist}.
\end{proof}
\end{Proposition}

Note that for a low herding intensity there exist one maximum at the average opinion of 0, while for a high herding intensity two symmetrical maxima emerge. 

\subsubsection{Large market approximation}

Although the stationary distribution in Equation (\ref{stationary_dist}) is analytically exact,  the calculation is demanding due to the binominal coefficient. To illustrate the usefulness of Theorem \ref{Diff_approx} we approximate $\overline{V}_t^n$ using the large market limit in the next Proposition.
 
\begin{Proposition}[Large Market Approximation]\label{DIff_example_lux1}\ \\

If the distribution of the initial average mood $F_{\overline{M}_0}^n$ has for $n \to \infty$ a limit $F_{\overline{M}_0}$, then 
\begin{equation}
(\overline{V}_t^n)_{t \in [0, \infty)} \xrightarrow{\mathcal{L}}  (\overline{V}_t)_{t \in [0, \infty)} \ in \ D_{[-1,1]}[0,\infty),
\end{equation}
where $(\overline{V}_t)_{t \in [0, \infty)} $ is the solution of the ODE
\begin{equation}\label{Ode_V_t}
d\overline{V}_t =  2 \beta\left[ \tanh(\gamma \overline{V}_t) - \overline{V}_t \right]\cosh(\gamma \overline{V}_t)dt, \ \overline{V}_0 = \overline{\theta},
\end{equation}
with $\overline{\theta} \sim F_{\overline{M}_0}$.

\begin{proof}
See Appendix \ref{proof_DIff_example_lux1}
\end{proof}
\end{Proposition}

\begin{figure}[h!]
{\centering
\includegraphics[width=1\textwidth, height=4.5cm]{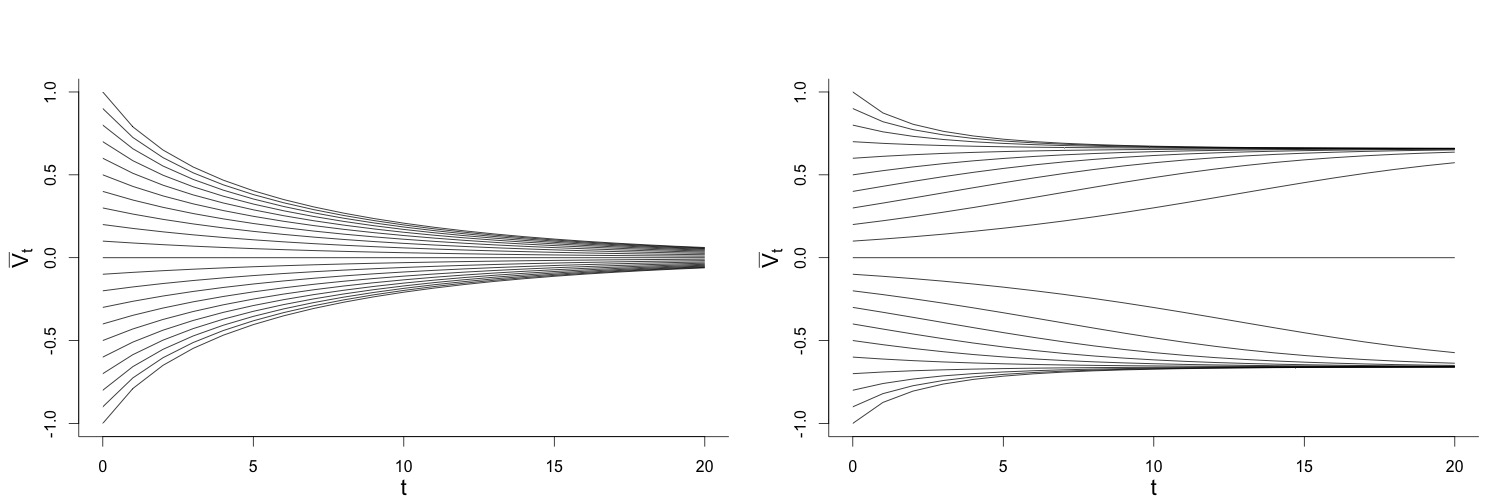}
\captionof{figure}{$\overline{V}_t$ for $\gamma$ = 0.8, 1.2 and $\beta$ = 0.3 with different $\overline{V}_0$}\label{P_st_sol_ode}
}
\end{figure}

Note, that the large market limit has different properties than the original market, i.e. no state transition and one, respectively two, absorbing states. 
Below we show the solution of the above ODE for different initial values $\overline{\theta}$. The left graph of Figure \ref{P_st_sol_ode} shows  $\overline{V}_t$ for $\gamma = 0.8$ and $\beta = 0.3$. Independent of the initial value $\overline{V}_0$, $\overline{V}_t$ converges monotone to 0 for $t \to \infty$. When $\gamma >1$,as illustrated in the right graph, the solution of Equation (\ref{Ode_V_t}) has three limits for $t \to \infty$ depending on the initial value $\overline{\theta}$. For $\overline{\theta} = 0$, $\overline{V}_t$ is constantly 0. For $\overline{\theta}  \underset{(<)}{>} 0$, $\overline{V}_t$ converges monotone to the positive (negative) solution of the Equation $y = \tanh(\gamma y)$.

Before we link the endogenous model to the price process in the next sub-section we show an application of Proposition \ref{prop_rate_of_conv}, i.e. we determine the rate of convergence.

\begin{Remark}[Rate of convergence]\label{remark_convergence_rate}\ \\
If $F_{\overline{M}_0}^n = F_{\overline{M}_0}$, then
\begin{equation}\label{inequality_convergence}
\sup_{s \leq t} |\overline{V}_s^n-\overline{V}_s| \leq 2 \sqrt{\frac{\beta}{n}} \sup_{s \leq t} |Y_s| \leq \frac{C_{\beta,\gamma}}{\sqrt{n}}, \ \forall t \geq 0,
\end{equation}

where $Y_t$ is the unique solution of the SDE

\begin{equation}\label{SDE_Y}
dY_t = \sqrt{(1-Y_t \tanh(\gamma Y_t))\cosh(\gamma Y_t)} dB_t, \ Y_0 = 0,
\end{equation}

and $C_{\beta,\gamma}$ is a constant depending on $\beta$ and $\gamma$ via $C_{\beta,\gamma}=2\sqrt{\beta}| y^*| $, where $y^*$ is the solution of $1-y\tanh(\gamma y) = 0$.

\begin{proof}
Let $a_n := \sqrt{\frac{n}{4 \beta}}$. Moreover define $\hat{c}(y)^2 := (1-y\tanh(\gamma y))\cosh(\gamma y)$ and $\hat{b} := \overline{b}$. Now, the first inequality in  (\ref{inequality_convergence}) simply follows from application of Proposition \ref{prop_rate_of_conv}.\\
Since $\forall \gamma > 0$, $\hat{c}(y)$ is symmetric to zero from which it falls monotonously to 0, the related process $Y_t$ has less variance the more its distance to zero. When the process $Y_t$ reaches $y^*$, where $y^*$ is the solution of $\hat{c}(y^*) = 0$, the diffusion coefficient is zero and $Y_t$ stays constant. Hence $\sup_{s \leq t} |Y_s| \leq |y^*|, \ \forall t \geq 0$, which concludes the second inequality.
\end{proof}

\end{Remark}

\newpage

\subsection{Price dynamics}

In this sub-section we extend the endogenous model of the previous sub-section with a link to an asset price. Therefore we characterize the impact of agents opinion on the asset price and vice versa.   
Moreover we introduce an additional group of traders called fundamentalists.\footnote{Note that, while in the model of Lux \cite{lux1995} fundamentalists are required to instantly match supply and demand, we could forgo, since in our model orders arrive asynchronously.} Last are characterized by basing their behavior on the difference between actual price and a \textit{fundamental value} $F \in \mathbb{R}$. In particular, when the price is below (above) $F$, they consider the asset cheap (expensive) and want to buy (sell).
We assume the fundamentalists are homogeneous, viz. $F$ is common, and the fundamental value is time-invariant.
Although we are mostly consistent with Lux \cite{lux1995} in setting our assumptions, we additionally introduce random signals reflected in the agents trading behavior.\\
Let $k_n \in \mathbb{N}$ denote the number of fundamentalists and $\phi_n = \frac{k_n}{n} $ the proportion of fundamentalists within all agents.
Hence $\mathbb{A}_n = \{1, ... , n \}, \ \ n-k_n \in 2 \mathbb{N}$ is the set of all agents participating in the market.

Compared to the previous sub-section we extent the state space to be $S = \{-1,1,2 \}$, where $s_3 = 2$ denotes a fundamentalistic agent. 
The market character is then given by $M_k = (M_k^1, M_k^2, M_k^3)$, where $M_k^1$ and $M_k^2$ is defined in the previous subsection (Equations (\ref{M_k^1_Lux}) and (\ref{M_k^2_Lux})) and
\begin{equation}
M_k^3 = \frac{1}{n} \sum_{a \in \mathbb{A}_n} \mathds{1}_{\{2\}}(x_k^a), \ k \geq 0.
\end{equation}

In line with Lux we define the \textit{average opinion of noise traders} by

\begin{equation}
\overline{M}_k := \left( \frac{1}{1-\phi_n} \right) (M_k^2-M_k^1), \ k \geq 0.
\end{equation}

Moreover let $\widetilde{F}_{x_0}^n$ denote the initial distribution of states and $\widetilde{F}_{M_0}^n$ the resulting initial distribution of $M_0$.
We assume fundamentalists weight their demand according to a to $F$ symmetric function $w_2: \R \to \R$ depending on the current price and $F$. On the other hand optimists buy a fixed amount of shares ($w_1$) while pessimists want to sell the same amount. 
We scale the demand by $1/\sqrt{n}$ and add a random signal $\xi_k$, which is assumed to be i.i.d. with $\E[\xi_1]=0$ and $\Var[\xi_1]<\infty$. Let $\widetilde{F}_{P_0}$ be the distribution of the starting price, which for simplicity is assumed independent of $n$.

\begin{Definition}[Excess demand function]\ \\
In summary we assume the following excess demand function
\begin{equation}
e_a^n(P_{k-1},\xi_k)=\begin{cases}
n^{-1/2} w_1 x_k^a + \xi_k, & x^a_k \in \{-1,1 \}\\
n^{-1/2} w_2(F,P_{k-1}) + \xi_k,  & x^a_k = 2,
\end{cases}
\end{equation}

\end{Definition}

Next we define the pricing rule, that is how supply or demands of agents impact the stock price. 

\begin{Definition}[Pricing rule]\ \\
We assume the pricing rule is given by
\begin{equation}
r_n(q,x) = x + \frac{\alpha}{\sqrt{n}} q
\end{equation}
\end{Definition}

In Lux's model the price dynamics are an equilibrium result of matching supply and demand of all participating agents.\footnote{Since all agents are considered by the supply and demand matching, a transition to an infinite big market (viz. $n \to \infty$) is not readily possible within Lux's framework.} Since we assumed that at a specific point in time only one agent is trading (see Assumption \ref{assumption_k-th_action}) and therefore impacts the price solely, the pricing rule in our model largely deviates from Lux by construction. However, we can use the trading intensity functions $\lambda_a$ in order to "scale" the number of trades and make the models comparable within a fixed time interval. More precisely, following the homogeneity assumption related to our agents, we set 

\begin{equation}
\lambda_a(x,v) = \bar{\lambda} \in \mathbb{R}^+
\end{equation}

Note that for $\bar{\lambda} = n$ the expected net excess demand of all agents for a fixed time interval is the same within our and Lux's model.

Now, also to incorporate a feedback effect from the price to agents behavior we extent the transition probability of the mood-based traders. Therefore we let the transition probability not only depend on the overall mood, but also on the expected price dynamics.
\begin{Definition}[Transition probabilities]\ \\
The transition probability to switch the state from $-1$ to $1$ is defined to be
\begin{equation}\label{Transition probability Lux1}
\Pi^{1,2}_n(P_{k-1},M_{k-1}) = \beta e^{\gamma_1 \widehat{z}_n(P_{k-1},M_{k-1}) + \gamma_2 \overline{M}_{k-1}},
\end{equation}

where
\begin{equation}
\widehat{z}_n(P_{k-1},M_{k-1}) := \overline{\lambda} \left( \frac{k_n}{n} w_2 (F,P_{k-1}) + \left(1-\frac{k_n}{n}\right) w_1 \overline{M}_{k-1}\right).
\end{equation}

The transition probability to switch the state from $1$ to $-1$ is defined as
\begin{equation}\label{Transition probability Lux2}
\Pi^{2,1}_n(P_{k-1},M_{k-1}) = \beta e^{-\gamma_1 \widehat{z}_n(P_{k-1},M_{k-1}) - \gamma_2 \overline{M}_{k-1}},
\end{equation}
where $\gamma_1, \gamma_2 >0$. Moreover we assume that fundamentalist can not become optimists/ pessimists and vice versa.\footnote{Note, that the transition probabilities in Equations (\ref{Transition probability Lux1}) and (\ref{Transition probability Lux2}) are not per se well defined. Instead of capping the probabilities at one we rather use the function $w_2$ in order to control the impact of large prices $P_{k-1}$} Hence the transition matrix is given as
\begin{equation}
\Pi_n(P_{k-1},M_{k-1}) = \begin{pmatrix}  1-\Pi^{1,2}_n & \Pi^{1,2}_n &  0 \\ \Pi^{2,1}_n & 1-\Pi^{2,1}_n  &  0 \\ 0 & 0 & 1   \\ \end{pmatrix} (P_{k-1},M_{k-1}).
\end{equation}
\end{Definition}

While $\widehat{z}_n$ measures the expected price dynamics, $\gamma_1$ measures the intensity of the price feedback on agents behavior. On the other hand $\gamma_2$ describes the herding intensity analogous to the previous chapter.\\ 
Moreover since $M_k^3 = \phi_n$ and $M_k^2 = (1-\phi_n) - M_k^1$ the market character is uniquely defined by the average opinion of the noise traders, i.e.
 
\begin{equation}\label{M_k as overline M_k}
M_k = \left( \frac{(1-\phi_n)(1-\overline{M}_k)}{2}, \frac{(1-\phi_n)(1+\overline{M}_k)}{2}, \phi_n \right)
\end{equation}

Analogous to Equation (\ref{def_average_mood_process}) and (\ref{M_k as overline M_k}) the average opinion is derived as
\begin{equation}
\overline{V}_t^n = \sum^\infty_{k=0} \overline{M}_k \mathds{1}_{[T_k,T_{k+1})}(t) = \frac{1}{1-\phi} \left(V^{n,2}_t - V^{n,1}_t \right), \ t \geq 0.
\end{equation}

Now, to determine the behavior of the exemplary model we again leverage from the results presented in Lux \cite{lux1995}. Although our model is different by construction, the key factors like net-excess demand, weighting of fundamentalists and mood traders, etc. are comparable.
In the next remark we state the behavior of our model, which is valid not only for the price process $X_t^n$ and the average opinion index $\overline{V}_t^n$, but also for the underlying markov chains $(P_k)_{k \geq 0}$ and $(\overline{M}_k)_{k \geq 0}$. 
 
\begin{Remark}[Market behavior]\label{approx_beh_ext_lux}\ \\
\begin{enumerate}
\item For a high herding intensity $\gamma_2$, there exist two equilibria $E_+ = (\widetilde{v}_+,x_+)$ and $E_- = (\widetilde{v}_-,x_-)$, where $\widetilde{v}_+ = -\widetilde{v}_-$ and $F-x_- = x_+ - F$.
\item For a small herding intensity $\gamma_2$, there is one unique equilibrium $E_0=(0,F)$. ,
\begin{enumerate}
\item If the intensity of price feedback $\gamma_1$ is low, then $E_0$ is stable. 
\item If the intensity of price feedback $\gamma_1$ is large, $E_0$ is unstable and there occur periodic cycles.
\end{enumerate}
  
\end{enumerate}

\end{Remark}

\begin{figure}[h!]
{\centering
\includegraphics[width=1\textwidth, height=6cm]{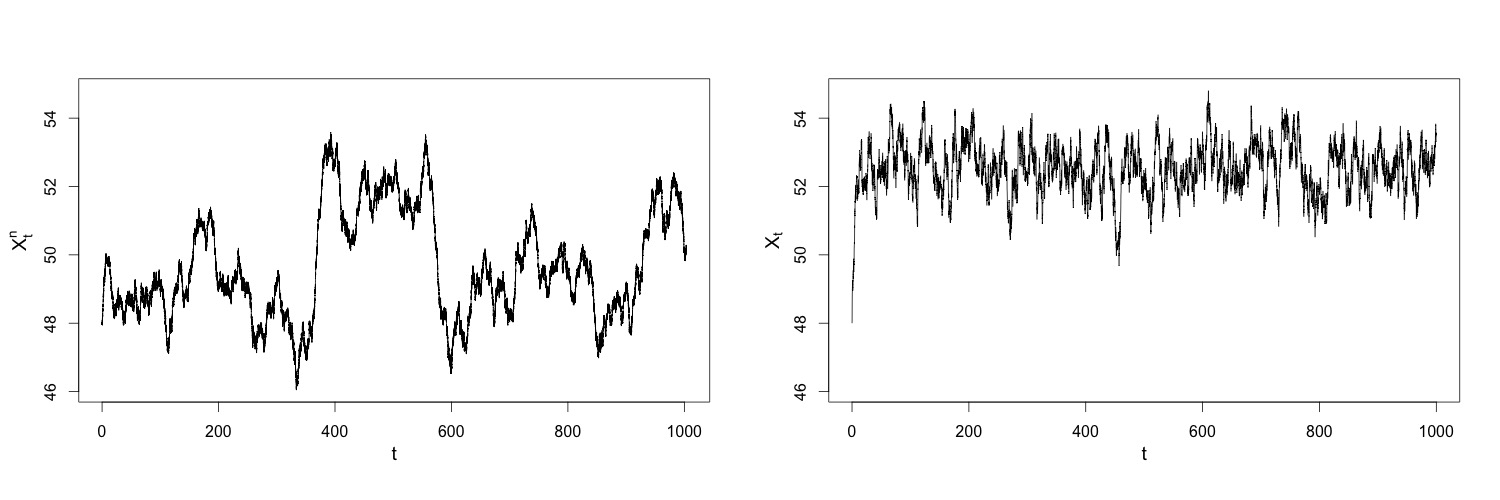}
\captionof{figure}{$X^n_t$ and $X_t$ for $\gamma_1$ = 0.2, $\gamma_2$ = 1.2, $w_2 = F - x$}\label{Fig_Xn_X_1}
\medskip
\includegraphics[width=1\textwidth, height=6cm]{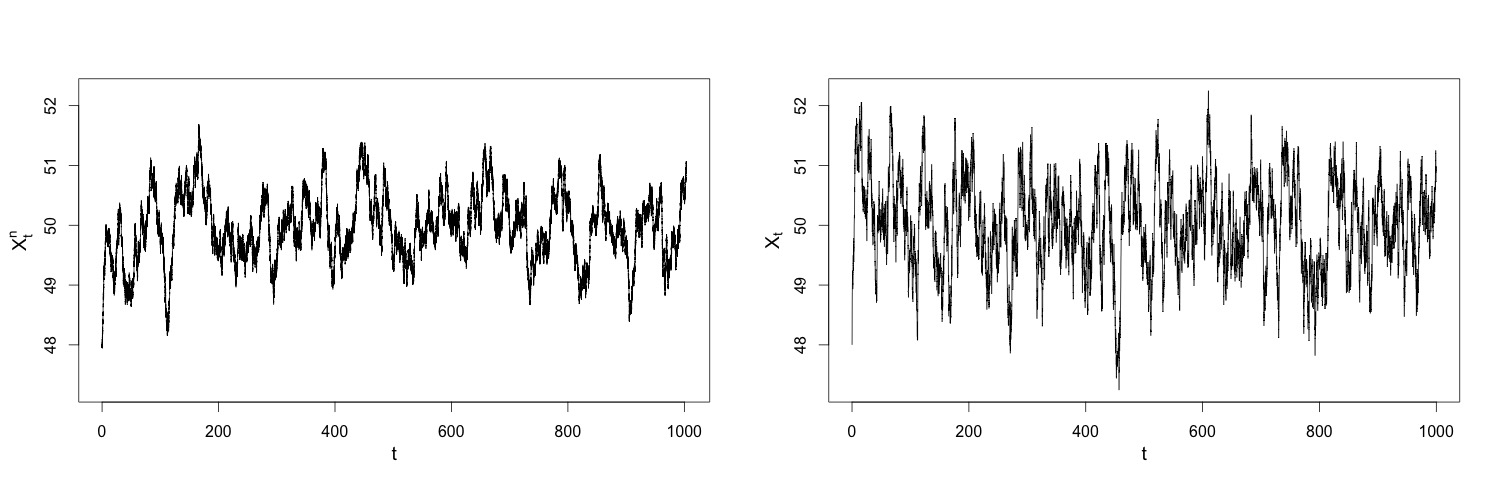}
\captionof{figure}{$X^n_t$ and $X_t$ for $\gamma_1$ = 0.2, $\gamma_2$ = 0.8, $w_2 = F - x$}\label{Fig_Xn_X_2}
\medskip
\includegraphics[width=1\textwidth, height=6cm]{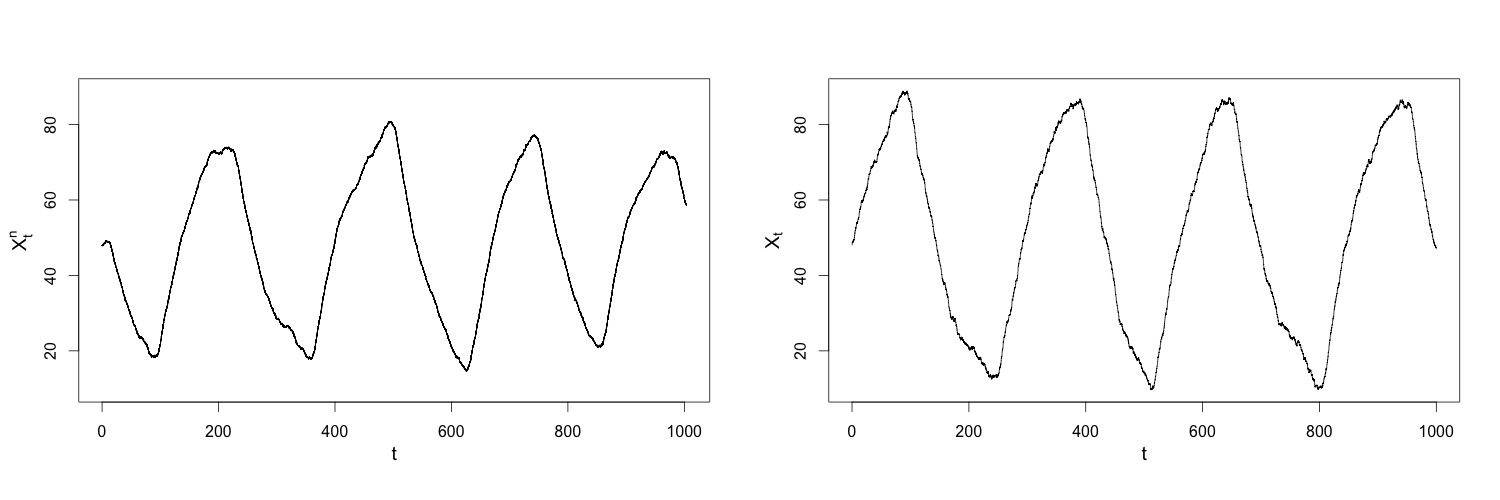}
\captionof{figure}{$X^n_t$ and $X_t$ for $\gamma_1$ = 1.2, $\gamma_2$ = 0.8, $w_2 = 0.05*(F - x)$}\label{Fig_Xn_X_3}
}
\end{figure}

Although we refrain from a more detailed description of the market behavior within this paper, we illustrate Remark \ref{approx_beh_ext_lux} by showing trajectories of $X_t^n$ for each case in Figure \ref{Fig_Xn_X_1} - Figure \ref{Fig_Xn_X_3}.\\  

Also for the extended example we state the large limit approximation, which is in contrast to the pure endogenous large market dynamics a diffusion price process.

\begin{Proposition}[Large market approximation]\label{Diff_ex2}\ \\

If $\phi_n \xrightarrow{n \to \infty} \phi$ and $\overline{V}_0^n \xrightarrow{\mathcal{L}} \overline{\theta}$, then

\begin{equation}
(X_t^n,V_t^n)_{t \in [0, \infty)} \xrightarrow{\mathcal{L}}  (X_t,V_t)_{t \in [0, \infty)} \ in \ D_{\mathbb{R} \times [0,1]^2}[0,\infty),
\end{equation}
where $(X_t,V_t)_{t \in [0, \infty)}$ is the unique strong solution of the SDE
\begin{equation}\label{Eq_SDE_extended}
\begin{cases}
dX_t = \alpha z(X_t,\overline{V}_t) dt + \alpha (\bar{\lambda} \Var[\xi_1])^{1/2} dB_t, &X_0 = \zeta\\
d\overline{V}_t = 2 \beta \left[ \tanh ( \gamma_1 z(X_t,\overline{V}_t) + \gamma_2 \overline{V}_t) - \overline{V}_t \right] \cosh(\gamma_1 z(X_t,\overline{V}_t) + \gamma_2 \overline{V}_t)dt, & \overline{V}_0 = \overline{\theta}\\
dV_t^3 = 0, & V^3_0 = \phi,
\end{cases}
\end{equation}
 
where $\overline{V}_t = \frac{1}{1-\phi} (V^2_t - V^1_t)$, $z(x,v) := \overline{\lambda}[\phi w_2(F,x) + (1-\phi) w_1 v]$ and $\zeta \sim \widetilde{F}_{P_0}$.  
 
\begin{proof}
See Appendix \ref{proof_DIff_ex2}
\end{proof}
\end{Proposition}

\newpage
 
\newpage

In Figure \ref{Fig_Xn_X_1} - Figure \ref{Fig_Xn_X_3} we show a trajectory of the solution of Equation (\ref{Eq_SDE_extended}) with $\beta = 0.12$, $\phi = 0.2$, $w_1=1$, $P_0 = 48$, $F=50$ and $\sqrt{\Var[\xi_1]}=0.2$. We compare $X_t^n$ for $n=100$ and $X_t$ in 1000 time units. As such we compare a medium large market with the infinite large market, which is used to approximate. 
In Figure \ref{Fig_Xn_X_1} we set $\gamma_1 = 0.2$, $\gamma_2 = 1.2$ and $w_2 = F-x$.  In line with Remark \ref{approx_beh_ext_lux} 1. we see a regime switch between two equilibria which are symmetric to the fundamental value for $X_t^n$. The transition also holds true for the diffusion process $X_t$ when $\sqrt{\Var[\xi_1]}>0$, however on a much larger scale than shown in Figure \ref{Fig_Xn_X_1}. For completeness we illustrate the transition of $X_t$ on a 100 times larger scale in the Appendix with Figure \ref{Extension example_figure}.
If we now reduce the herding intensity $\gamma_2$ to 0.8, as illustrated in Figure \ref{Fig_Xn_X_2}, $X_t^n$ as well as $X_t$ have an equilibrium at the same point, namely the fundamental value $F$, which is in accordance with Remark \ref{approx_beh_ext_lux} 2.(a). Nevertheless the diffusion process shows a higher stochastic noise.
In Figure \ref{Fig_Xn_X_3} we illustrate $X_t^n$ and $X_t$  in the case that the intensity of price feedback is high while the influence of fundamentalists is low. Therefore we set $\gamma_1$ to 1.2 and $w_2$ to $0.05*(F-x)$. Independent of the initial distribution of optimists and pessimists, $X_t^n$ and $X_t$ are then oscillating around $F$ with the same scale, although the amplitude of the diffusion process is slightly higher.
In summary, the diffusion process $X_t$ shows the same characteristics as $X_t^n$ and is well suited to be used as large market approximation to examine those. Nevertheless, as apparent from Figure \ref{Fig_Xn_X_1} the scaling in which the characteristic is displayed, might be different.

\section{Conclusion and Outlook}
We provided a microscopic model that can be used not only to study markovian socio-economic dynamics and related price processes but also diffusion price processes resulting as a large market limit. As such we extend the model of Pakkanen \cite{pak2010} with state space based socio-economic behavior. After we connected the endogenous dynamics with the price process within a finite Markovian environment, we establish the circumstances under which the finite market can be approximated by an (It{\^o}-) diffusion process. Although our model is rather general compared to other models available (see Introduction), we still make the strong assumption of Markov property. This seems unrealistic, especially in light of financial markets, where agents may follow trends and adjust their behavior from experience. Already relaxing the Markov assumption a bit leads to a way more complex theory behind. Assume, for example, that agents consider a fixed time of experience in their behavior. As for future research one might show that with this assumption the finite market can be seen as a micro-foundation for price processes that are a solution of an stochastic differential delay equation, whose use is rather new in mathematical finance (see e.g.  Arriojas, Mercedes, et al. \cite{arr2008}).\\
To demonstrate the applicability of a seperation of behavior and price process we embedded the model of Lux \cite{lux1995} within our framework. We confirmed, using similar assumption as made in Lux \cite{lux1995}, that herding behavior can induce phase transitions and oscillations in the finite-market price process. Introducing random signals, which influence the agents excess demand function, we were able to extent the result to diffusion price processes, which are the result of a large market limit. Hence we provided an agent behavior based explanation for intrinsic price cycles often seen in financial markets.

\newpage
\section{Appendix}

\subsection{Additional Figures}\label{additional_figures}

\begin{figure}[h!]
{\centering
\includegraphics[width=1\textwidth, height=4.5cm]{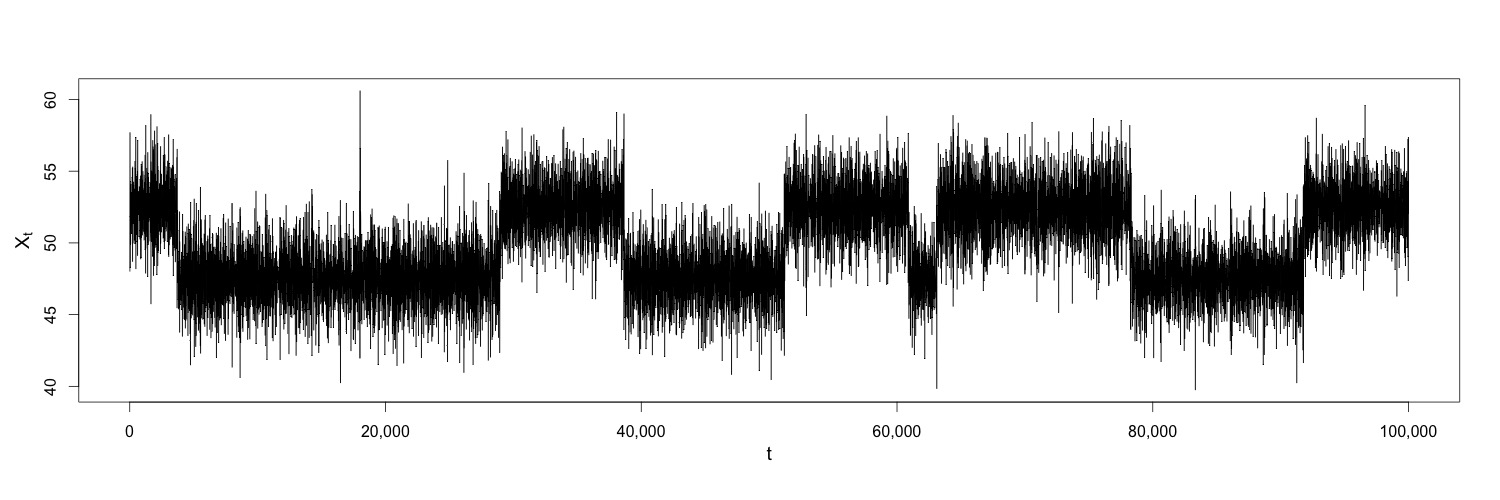}
\captionof{figure}{$X_t$ for $\gamma_1$ = 0.2, $\gamma_2$ = 1.2, $w_2 = F - x$}\label{Extension example_figure}
}
\end{figure} 


In the following, to study the dynamics behind the SDE presented in Equation \ref{Diff_ex2}, we look at the large market dynamics without random signals. Equation \ref{Diff_ex2} then simplifies to the ODE

\begin{equation}\label{Eq_ODE_extended}
\begin{cases}
dX_t = \alpha z(X_t,\overline{V}_t) dt, &X_0 = \zeta\\
d\overline{V}_t = 2 \beta \left[ \tanh ( \gamma_1 z(X_t,\overline{V}_t) + \gamma_2 \overline{V}_t) - \overline{V}_t \right] \cosh(\gamma_1 z(X_t,\overline{V}_t) + \gamma_2 \overline{V}_t)dt, & \overline{V}_0 = \overline{\theta}\\
dV_t^3 = 0, & V^3_0 = \phi,
\end{cases}
\end{equation}
 
where $\overline{V}_t = \frac{1}{1-\phi} (V^2_t - V^1_t)$, $z(x,v) := \overline{\lambda}[\phi w_2 (F,x) + (1-\phi) w_1 v]$ and $\zeta \sim \widetilde{F}_{P_0}$. \\ 
 
\begin{figure}[h!]
{\centering
\includegraphics[width=1\textwidth, height=6cm]{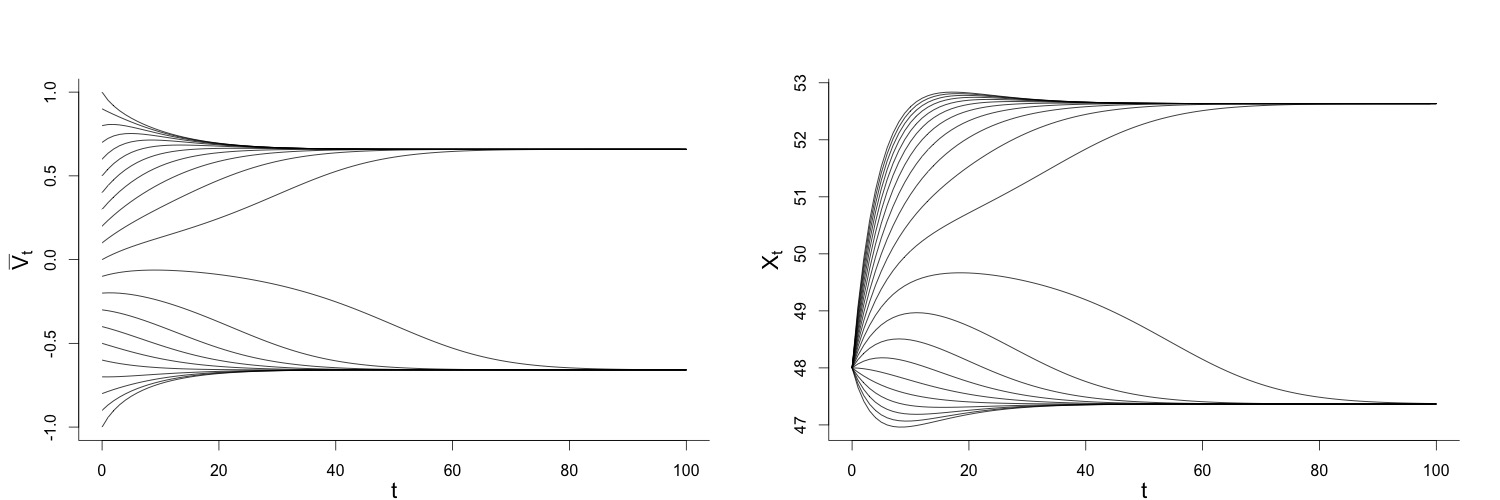}
\captionof{figure}{$\overline{V}_t$ and $X_t$ for $\gamma_1$ = 0.2, $\gamma_2$ = 1.2, $w_2 = F-x$}\label{Fig_ex_lux_ode_1}
\medskip
\includegraphics[width=1\textwidth, height=6cm]{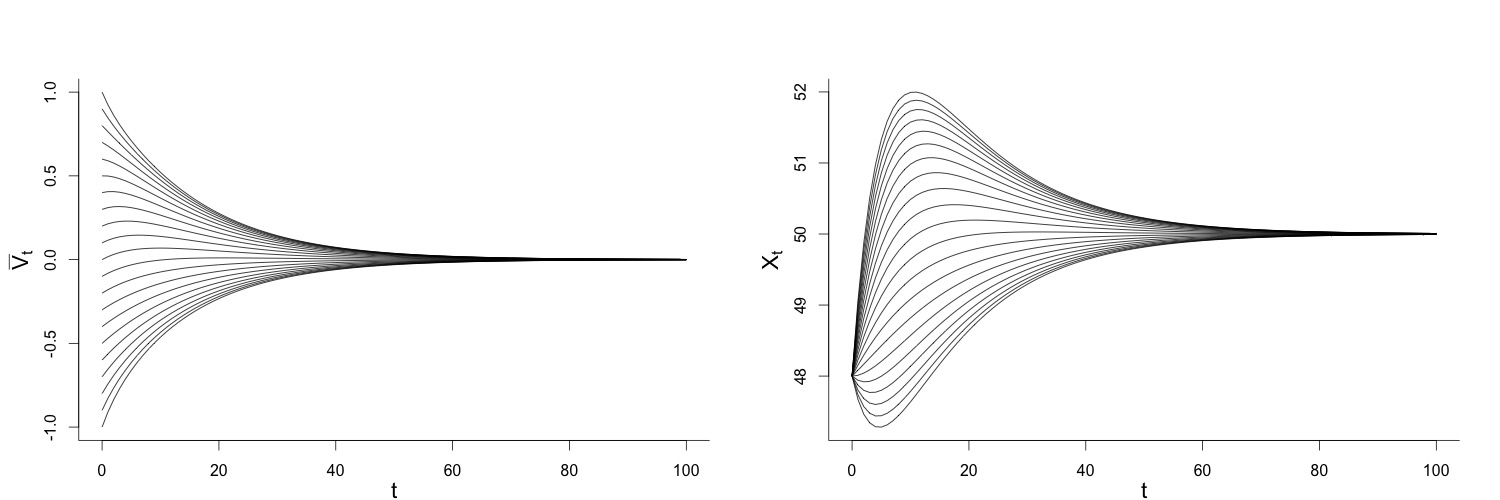}
\captionof{figure}{$\overline{V}_t$ and $X_t$ for $\gamma_1$ = 0.2, $\gamma_2$ = 0.8, $w_2=F-x$}\label{Fig_ex_lux_ode_2}
\medskip
\includegraphics[width=1\textwidth, height=6cm]{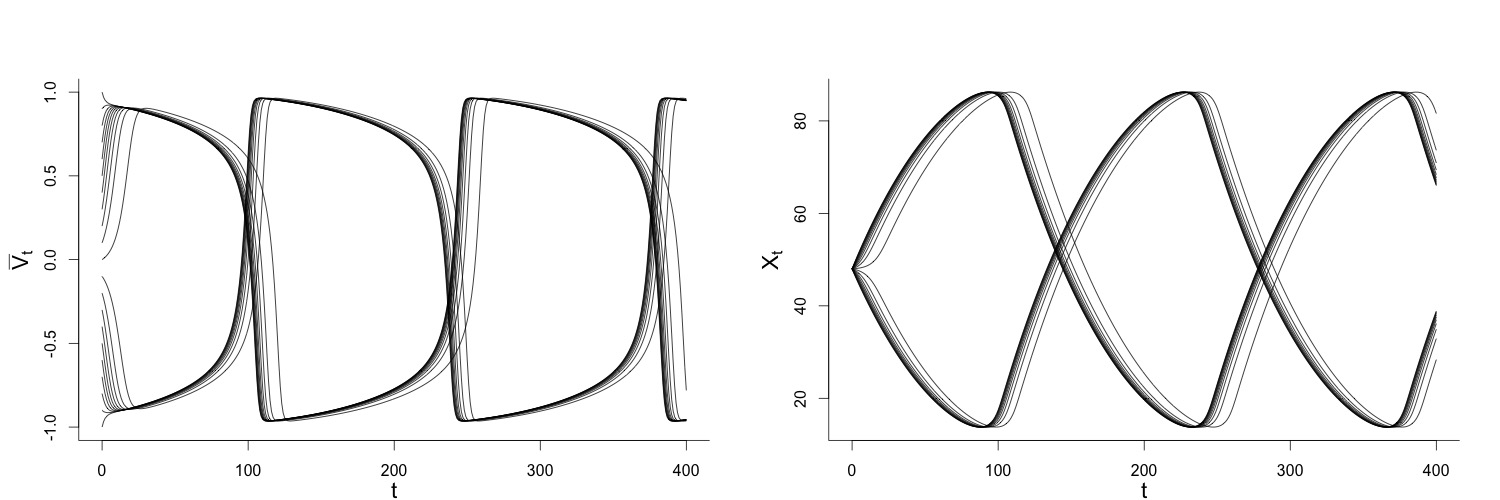}
\captionof{figure}{$\overline{V}_t$ and $X_t$ for $\gamma_1$ = 1.2, $\gamma_2$ = 0.8, $w_2= 0.05*F-x$}\label{Fig_ex_lux_ode_3}
}
\end{figure}

In Figure \ref{Fig_ex_lux_ode_1} - Figure \ref{Fig_ex_lux_ode_3} we show several solutions of Equation (\ref{Eq_ODE_extended}) for the different initial values $\overline{\theta}$ for the same set of parameters as in Figure \ref{Fig_Xn_X_1} - \ref{Fig_Xn_X_3}, that is $\beta = 0.12$, $\phi = 0.2$, $w_1=1$ and $P_0 = 48$. Figure \ref{Fig_ex_lux_ode_1} shows $\overline{V}_t$ and $X_t$ for 100 time units with the same setting as used in Figure \ref{Fig_Xn_X_1}, i.e. $\gamma_1 = 0.2$, $\gamma_2 = 1.2$ and $w_1 = 1$. Depending on the initial value $\overline{\theta}$, $\overline{V}_t$ and $X_t$  converge monotonously to one of two constants, which are symmetrical to 0, respectively F. If we now reduce the herding intensity $\gamma_2$ to 0.8, as illustrated in Figure \ref{Fig_ex_lux_ode_2}, $\overline{V}_t$ and $X_t$  converge monotonously to 0, respectively F, independent of the initial value. In \ref{Fig_ex_lux_ode_3} we illustrate $\overline{V}_t$ and $X_t$  in the case that the intensity of price feedback is high while the influence of fundamentalists is low. Therefore we set $\gamma_1$ to 1.2 and $w_2$ to $0.05*(F-x)$. Independent of the initial distribution of optimists and pessimists, $\overline{V}_t$ and $X_t$ are then oscillating around 0 , respectively $F$.

As observable in Figure \ref{Fig_ex_lux_ode_1} and Figure \ref{Fig_ex_lux_ode_2} the solution of Equation (\ref{Eq_ODE_extended}) may converge to constants $x$, $\overline{v}$. In order to specify the constants we require $z(X_t,\overline{V}_t) = 0$ and $d\overline{V}_t = 0 $, which is equivalent to 
\begin{equation}\label{Eq_ODE_extended2}
\begin{cases}
\phi w_2(F,x) + (1-\phi) w_1 \overline{v} = 0\\
\left[\tanh (\gamma_2 \overline{v}) - \overline{v} \right] \cosh(\gamma_2 \overline{v}) = 0
\end{cases}
\end{equation}
Hence with the parameters above, $x= 50 \pm 4\overline{v}$, where $\overline{v}$ is the solution of $y = \tanh(\gamma_2 y)$, if $\gamma_2 >1$ and 0 otherwise.

\newpage
 
\newpage

\subsection{Proof of Theorem \ref{Diff_approx}}\label{proof_diff_approx1}
In order to determine the diffusional limit we want to apply Theorem IX. 4.21 of Jacod and Shiryaev \cite{jac2003}. 
Let us first note that the problem is well defined as the local Lipschitz and linear growth conditions set related to $z$, $b$,$\sigma$ and $c$ imply the existence of a unique solution of Equation (\ref{SDE1}) by Theorem III.2.32 of Jacod and Shiryaev \cite{jac2003}. Moreover IX.4.3. (ii) Jacod and Shiryaev \cite{jac2003} follows from Theorem 21.10 Kallenberg \cite{kal2002}.\\
In order to improve the readability we capture the change of the price by 
\begin{equation}
\overline{r}_n(q,x) = r_n(q,x) - x
\end{equation}

\subsubsection*{Hypothesis (i)}
To show Theorem IX. 4.21 (i) of Jacod and Shiryaev \cite{jac2003} we have to calculate the first and second moments of the rate kernel $K_n$ defined in Equation (\ref{rate_kernel}).
We start with the first moment related to the price component. Applying Equations (\ref{prob_trade}),(\ref{rate_kernel}), (\ref{pricing_rule}) and the desintegration theorem stated in Theorem 6.4 of Kallenberg \cite{kal2002} leads to the following representation.

\begin{equation}
\begin{split}
\int K_n(x,v,dy,dw)y &= \nu_{\mathbb{A}_n}(x,v) \int k_n(x,v,dy,dw)y \\
&=\nu_{\mathbb{A}_n}(x,v) \E [P_1-P_0|P_0=x,M_0=v]\\
&= \sum_{a=1}^n \lambda_a(x,v) \E \left[\overline{r}_n(e_a^n(x,v,\xi_1),x)\right]\\
&= \alpha z_n(x,v) + \sum_{a=1}^n \lambda_a(x,v) \E [u_n(e_a^n(x,v,\xi_1),x)]
\end{split}
\end{equation}

Moreover by Assumption \ref{assumption_excess_demands} for any $\delta >0$ and $|(x,v)| < \delta$ we have
\begin{equation}
\begin{split}
&\left| \sum_{a=1}^n \lambda_a(x,v) \E [u_n(e_a^n(x,v,\xi_1),x)] \right|\\
&\leq n C_n^{\delta} \sup_{|(x,v)| < \delta} \left| \frac{\lambda_{\mathbb{A}_n}(x,v)}{n} \right|  \sup_{\substack{|(x,v)| < \delta, \\ \ a \in \mathbb{A}_n}} \E [|e_a^n(x,v,\xi_1)|],
\end{split}
\end{equation}
where $C^{\delta}_n = o(n^{-1})$ and the other terms are bounded as assumed in Assumption \ref{assumption_explosion}. Thus the u.o.c. convergence to zero when $n \to \infty$.

Using the desintegration theorem again, Equation (\ref{rate_kernel}) and the dynamics stated in Equations (\ref{dynamics_mood2}) - (\ref{dynamics_mood3}) together with representation (\ref{b_n^i}) we get the first moment related to the occupancy measure of the single states $s_i$ via

\begin{equation}
\begin{split}
&\int K_n(x,v,dy,dw)w_i = \nu_{\mathbb{A}_n}(x,v) \int k_n(x,v,dy,dw)w_i \\
&= \nu_{\mathbb{A}_n}(x,v) \E [M_1^i - M_0^i|P_0 = x, M_0 = v]\\
&= n^{-d_1} \sum_{a=1}^n \mu_a(x,v) (\Pi^{i+}_{n,a}(x,v) - \Pi^{i-}_{n,a}(x,v) )\\
&= b_n^i(x,v).
\end{split}
\end{equation}

Since we assumed $b_n \to b$ when $n \to \infty$ in Theorem \ref{Diff_approx}, in summary we have
\begin{equation}
\int K_n(x,v,dy,dw)(y,w) \xrightarrow[n \to \infty]{u.o.c.} (\alpha z(x,v), b(x,v))
\end{equation}

Now let us consider the second moments.  
\begin{equation}
\begin{split}
&\int K_n(x,v,dy,dw)y^2 =  \nu_{\mathbb{A}_n}(x,v) \int k_n(x,v,dy,dw)y^2 \\
&= \nu_{\mathbb{A}_n}(x,v) \E [(P_1-P_0)^2|P_0 = x, M_0 = v]\\
&=\sum_{a=1}^n \lambda_a(x,v) \E [\overline{r}_n(e_a^n(x,v,\xi_1),x)^2]\\
&= \alpha^2 \sigma_n(x,v)^2 + \rho_n(x,v)
\end{split}
\end{equation}

where
\begin{equation}
\rho_n(x,v) := \sum_{a=1}^n \lambda_a (x,v) \E \left[2 \alpha n^{-d_2} e^a_n(x,v,\xi_1) u_n^a(e_n^a(x,v,\xi_1),x) + u_n ^a(e_n^a(x,v,\xi_1),x)^2\right].
\end{equation}

Using again Assumption \ref{assumption_excess_demands} we have $\forall \delta > 0$ and $|(x,v)| < \delta$
\begin{equation}
\begin{split}
|\rho_n(x,v)| &= \left| \sum_{a=1}^n \lambda_a (x,v) \E \left[ 2 \alpha n^{-d_2} e^a_n(x,v,\xi_1) u_n^a(e_n^a(x,v,\xi_1),x) + u_n ^a(e_n^a(x,v,\xi_1),x)^2\right] \right|\\
&\leq \sup_{\substack{|(x,v)| < \delta, \\ \ a \in \mathbb{A}_n}} \E [(2 \alpha n^{-d_2} C_n^{\delta} + (C_n^{\delta})^2) |e_n^a(x,v,\xi_1)|^2] \left| \sum_{a=1}^n \lambda_a (x,v) \right| \\
&\leq (2 \alpha n^{1-d_2} C_n^{\delta} + n(C_n^{\delta})^2) \sup_{\substack{|(x,v)| < \delta, \\ \ a \in \mathbb{A}_n}} \E[|e_n^a(x,v,\xi_1)|^2] \sup_{|(x,v)| < \delta} \left| \frac{\lambda_{\mathbb{A}_n}(x,v)}{n} \right|.
\end{split}
\end{equation}

Now, $\rho_n$ vanishes for $n \to \infty$ as $2 \alpha n^{1-d_2} C_n^{\delta} + n(C_n^{\delta})^2$ converges to zero and the other two terms are bounded by Assumption \ref{assumption_explosion}.

The second moments related to each state are given by

\begin{equation}
\begin{split}
&\int K_n(x,v,dy,dw)w_i^2 = \nu_{\mathbb{A}_n}(x,v) \E \left[ (M_1^i - M_0^i)^2|P_0 = x, M_0 = v \right]\\
&=n^{-2d_1}\sum_{a=1}^n \mu_{a}(x,v)(\Pi^{i+}_{n,a}(x,v) + \Pi^{i-}_{n,a}(x,v) )\\
&= c^i_n(x,v)^2
\end{split}
\end{equation}
which converges to $c^i(x,v)^2$ by Assumption \ref{var_states}.\\

As we assumed that the agents are not able to change their state and trade at the same time we have
\begin{equation}
\int K_n(x,v,dy,dw)y w_i = 0, \ \forall 1 \leq i \leq m.
\end{equation}

Moreover as only one agent at a time can change his state from $s_i$ to $s_j$ and as only the two respective state occupancy measures are affected by a change we have
\begin{equation}
\begin{split}
\int K_n(x,v,dy,dw)w_i w_j &= \nu_{\mathbb{A}_n}(x,v) \E \left[(M_1^i - M_0^i)(M_1^j - M_0^j) |P_0 = x, M_0 = v\right]\\
&= -n^{-(d_1+1)}\ \sum_{a=1}^n \mu_{a}(x,v)\left( v_i \Pi^{i,j}_{n,a}(x,v) + v_j \Pi^{j,i}_{n,a}(x,v) \right)\\
&= c^{i, j}_n(x,v)^2,
\end{split}
\end{equation}
which converges to $c$ for $n \to \infty$ as assumed in Theorem \ref{Diff_approx}.

Hence in summary we have 

\begin{equation}
\begin{split}
\int K_n(x,v,dy,dw) \begin{bmatrix} y^2 & yw   \\ wy   & w^2 \\ \end{bmatrix} &= \begin{bmatrix} \alpha^2 \sigma_n(x,v)^2 + \rho_n(x,v) & 0   \\ 0   & c_n(x,v)^2 \\ \end{bmatrix} \\
&\xrightarrow[n \to \infty]{u.o.c.} \begin{bmatrix} \alpha^2 v(x,v)^2  & 0   \\ 0   & c(x,v)^2\\ \end{bmatrix}
\end{split}
\end{equation}
where
\begin{equation}
\begin{bmatrix} y^2 & yw   \\ wy   & w^2 \\ \end{bmatrix} := \begin{bmatrix}
 y^2   & yw_1   & ... & yw_m\\
 w_1y & w_1^2 &  ... & w_1 w_m\\
 \vdots & \vdots & \ddots & \vdots \\
 w_m y & w_m w_1 & ... & w_m^2\\
\end{bmatrix}
\end{equation}

\subsubsection*{Hypothesis (ii)}
Next we show part (ii) of Theorem IX. 4.21 in Jacod and Shiryaev \cite{jac2003} i.e.

\begin{equation}
\sup_{|(x,v)|< \delta} \int K^n(x,v,dy,dw) |(y,w)|^2 \mathds{1}_{\{|(y,w)|>\epsilon\}} \xrightarrow{n \to \infty} 0 , \ \forall \epsilon > 0
\end{equation}

Therefore let $\epsilon_1 > 0$ and $\delta > 0$. Thus for $|(x,v)| < \delta$ we have
\begin{equation}
\begin{split}
&\sup_{|(x,v)|< \delta} \int K^n(x,v,dy,dw) y^2 \mathds{1}_{\{|y|>\epsilon_1\}}\\
&= \sup_{|(x,v)|< \delta} \sum_{a=1}^n \lambda_a(x,v) \E [\overline{r}_n(e_a^n(x,v,\xi_1),x)^2 \mathds{1}_{\{|\overline{r}_n(e_a^n(x,v,\xi_1),x)|> \epsilon_1 \}}]\\
&\leq \sup_{|(x,v)|< \delta} \left| \frac{\lambda_{\mathbb{A}_n}(x,v)}{n}\right| (\alpha n^{1/2 - d_2} + \sqrt{n} C_n^{\delta})^2 \sup_{\substack{|(x,v)| < \delta, \\ \ a \in \mathbb{A}_n}}  \E [|e_n^a(x,v,\xi_1)|^2\mathds{1}_{\{|e_a^n(x,v,\xi_1)|> \hat{\epsilon}_1^n \}}]\\
\end{split}
\end{equation}

where in the last inequality we used

\begin{equation}\label{ineq_r_n}
\begin{split}
\overline{r}_n(q,x) &= \alpha n^{-d_2} q + u_n^a(q,x)\\
&\leq \alpha n^{-d_2} |q| + \sup_{|x|< \delta} u_n^a(q,x)\\
&\leq (\alpha n^{-d_2} + C_n^{\delta} )|q|
\end{split}
\end{equation}

given by Assumption \ref{assumption_excess_demands} in the sense that 

\begin{equation}
\overline{r}_n(e_a^n(x,v,\xi_1),x)^2 \leq \frac{1}{n} (\alpha n^{1/2-d_2}+ \sqrt{n} C_n^\delta )^2 e_n^a(x,v,\xi_1)^2
\end{equation}
and as a result
\begin{equation}
|\overline{r}_n(e_a^n(x,v,\xi_1),x)| > \epsilon_1 \Leftrightarrow |e_n^a(x,v,\xi_1)| > \hat{\epsilon}_1^n
\end{equation}
with 
\begin{equation}
\hat{\epsilon}_1^n := \frac{\sqrt{n} \epsilon_1}{\alpha n^{1/2-d_2} + \sqrt{n} C_n^\delta}
\end{equation}

Now $\sup_{|(x,v)|< \delta} \left| \frac{\lambda_{\mathbb{A}_n}(x,v)}{n}\right| $ is bounded by Assumption \ref{assumption_explosion} and  $(\alpha n^{1/2 - d_2} + \sqrt{n} C_n^{\delta})^2$ converges to zero. Moreover  $C_n=o(n^{-1})$ and $\sup_{\substack{|(x,v)| < \delta, \\ \ a \in \mathbb{A}_n}}  \E [|e_n^a(x,v,\xi_1)|^2\mathds{1}_{\{|e_a^n(x,v,\xi_1)|> \hat{\epsilon}_1^n \}}]$ converges to zero by $\hat{\epsilon}_1^n  \xrightarrow{n \to \infty} \infty $ and uniform integrablity assumed in Assumption \ref{assumption_explosion}.
\\
Moreover let $\epsilon_2 > 0$ and $\delta > 0$. For $|(x,v)| < \delta$ we have
\begin{equation}
\begin{split}
&\sup_{|(x,v)|< \delta} \int K^n(x,v,dy,dw) w_i^2 \mathds{1}_{\{|w|>\epsilon_2\}}\\
&= \sup_{|(x,v)|< \delta} n^{-2d_1} \sum_{a=1}^n \mu_{a}(x,v)\left(\Pi^{i+}_{n,a}(x,v) + \Pi^{i-}_{n,a}(x,v) \right) \mathds{1}_{\{|w|>\epsilon_2\}}\\
&\leq 2 \sup_{|(x,v)|< \delta} n^{-2d_1} \sum_{a=1}^n \mu_{a}(x,v) \mathds{1}_{\{|w|>\epsilon_2\}}\\
&\leq 2 \sup_{|(x,v)|< \delta} \left| \frac{\mu_{\mathbb{A}_n}(x,v)}{n}\right| n^{1/2 - d_1} \mathds{1}_{\{|w|>\epsilon_2\}}\\
\end{split}
\end{equation}

and 

\begin{equation}
\begin{split}
&\sup_{|(x,v)|< \delta} \int K^n(x,v,dy,dw) |w_i w_j| \mathds{1}_{\{|w|>\epsilon_2\}}\\
&= \sup_{|(x,v)|< \delta} \left| -n^{-2d_1} \ \sum_{a=1}^n \mu_{a}(x,v)\left( n^{d_1-1} v_i \Pi^{i,j}_{n,a}(x,v) + n^{d_1-1} v_j \Pi^{j,i}_{n,a}(x,v) \right) \right| \mathds{1}_{\{|w|>\epsilon_2\}}\\
&\leq 2 \sup_{|(x,v)|< \delta} n^{-2d_1} \sum_{a=1}^n \mu_{a}(x,v) \mathds{1}_{\{|w|>\epsilon_2\}}\\
&\leq 2 \sup_{|(x,v)|< \delta} \left| \frac{\mu_{\mathbb{A}_n}(x,v)}{n}\right| n^{1/2 - d_1} \mathds{1}_{\{|w|>\epsilon_2\}}\\
\end{split}
\end{equation}

It is clear by Equation (\ref{dynamics_mood2}) and (\ref{dynamics_mood3}) that
\begin{equation}
R_n:= \operatorname{supp} (K_n(\cdot,v,dy,dw)(\cdot,w)) \subseteq \left[-n^{-d_1},n^{-d_1}\right]^m
\end{equation} 
Hence, since $R_n \cap \{ |w| \ > \epsilon_2 \} \xrightarrow{n \to \infty} \emptyset$ and $\sup_{|(x,v)|< \delta} \left| \frac{\mu_{\mathbb{A}_n}(x,v)}{n}\right|$ is bounded by Assumption \ref{assumption_explosion}, we see that
\begin{equation}
\sup_{|(x,v)|< \delta} \int K^n(x,v,dy,dw) |w|^2 \mathds{1}_{\{|w|>\epsilon_2\}} \xrightarrow{n \to \infty} 0.
\end{equation}

\subsection{Proof of Proposition \ref{prop_rate_of_conv}}\label{proof_prop_rate_of_conv}
In order to apply Corollary IX.4.28 of Jacod and Shiryaev \cite{jac2003} it is sufficient to show that $\forall \epsilon > 0$ and $\forall \delta > 0$ 
\begin{equation}
\sup_{|(x,v)|< \delta} a_n^2 \int K^n(x,v,dy,dw) |(y,w)|^2 \mathds{1}_{\{|y,w|>\epsilon / a_n \}} \xrightarrow{n \to \infty} 0.
\end{equation}

Let  $\epsilon_1 > 0$ and $\delta > 0$. Then for $|(x,v)| < \delta$ we have
\begin{equation}\label{eq_1_prop_rate_conv}
\begin{split}
& \sup_{|(x,v)|< \delta} a_n^2 \int K^n(x,v,dy,dw) |y|^2 \mathds{1}_{\{|y,w|>\epsilon / a_n \}}\\
&= \sup_{|(x,v)|< \delta} a_n^2 \sum_{a=1}^n \lambda_a(x,v) \E [\overline{r}_n(e_a^n(x,v,\xi_1),x)^2 \mathds{1}_{\{|\overline{r}_n(e_a^n(x,v,\xi_1),x)|> \epsilon_1/a_n \}}]\\
&\leq \sup_{|(x,v)|< \delta} \left| \frac{\lambda_{\mathbb{A}_n}(x,v)}{n}\right| (a_n(\alpha n^{1/2 - d_2} + \sqrt{n} C_n^{\delta}))^2 \sup_{\substack{|(x,v)| < \delta, \\ \ a \in \mathbb{A}_n}}  \E [|e_n^a(x,v,\xi_1)|^2\mathds{1}_{\{|e_a^n(x,v,\xi_1)|> \tilde{\epsilon}_1^n \}}]
\end{split}
\end{equation}

with $\tilde{\epsilon}_1^n := \frac{\sqrt{n} \epsilon_1}{a_n(\alpha n^{1/2-d_2} + \sqrt{n} C_n^\delta)}$. As we assumed that $\sqrt{n}a_n = O(n^{d_2} + (C_n^\delta)^{-1})$ we get $\tilde{\epsilon}_1^n \xrightarrow{n \to \infty} \infty$ and hence 

\begin{equation}
\{|e_a^n(x,v,\xi_1)|> \tilde{\epsilon}_1^n \} \xrightarrow{n \to \infty} \emptyset
\end{equation}
by the uniform integrablity assumed in Assumption \ref{assumption_trading_volume}. Moreover $\sup_{|(x,v)|< \delta} \left| \frac{\lambda_{\mathbb{A}_n}(x,v)}{n}\right| $ is bounded by Assumption \ref{assumption_explosion} which yield the convergence to zero of Equation (\ref{eq_1_prop_rate_conv}).

Moreover let  $\epsilon_2 > 0$ and $\delta > 0$. Then for $|(x,v)| < \delta$ we have
\begin{equation}
\begin{split}
& \sup_{|(x,v)|< \delta} a_n^2 \int K^n(x,v,dy,dw) |w_i|^2 \mathds{1}_{\{|w_i|>\epsilon_2 / a_n \}}\\
&= \sup_{|(x,v)|< \delta} a_n^2 \sum_{a=1}^n \mu_a(x,v) n^{-2 d_1}(\Pi_{n,a}^{i+} + \Pi_{n,a}^{i-}) \mathds{1}_{\{|n^{- d_1} (\Pi_{n,a}^{i+}-\Pi_{n,a}^{i-})|> \epsilon_2/a_n \}}\\
&\leq \sup_{|(x,v)|< \delta}  \left| \frac{\mu_{\mathbb{A}_n}(x,v)}{n}\right| 2 (a_n n^{1/2-d_1})^2 \mathds{1}_{\{a_n n^{- d_1} > \epsilon_2 \}}
\end{split}
\end{equation}

Now $\{a_n n^{- d_1} > \epsilon_2 \} \xrightarrow{n \to \infty} \emptyset$ and $a_n n^{1/2-d_1} < \infty$ since we assumed $a_n \sqrt{n} = O(n^{d_1})$. 

Moreover let  $\epsilon_3 > 0$ and $\delta > 0$. Then for $|(x,v)| < \delta$ we have
\begin{equation}
\begin{split}
& \sup_{|(x,v)|< \delta} a_n^2 \int K^n(x,v,dy,dw) |w_i w_j| \mathds{1}_{\{|w_i|>\epsilon_3 / a_n \}}\mathds{1}_{\{|w_j|>\epsilon_3 / a_n \}}\\
&= \sup_{|(x,v)|< \delta} a_n^2  \sum_{a=1}^n \mu_{a}(x,v) \left| -n^{-2d_1}  \left( n^{d_1-1} v_i \Pi^{i,j}_{n,a}(x,v) + n^{d_1-1} v_j \Pi^{j,i}_{n,a}(x,v) \right) \right| \\
&\cdot \mathds{1}_{\{|n^{- d_1} (\Pi_{n,a}^{i+}-\Pi_{n,a}^{i-})|> \epsilon_3/a_n \}} \mathds{1}_{\{|n^{- d_1} (\Pi_{n,a}^{j+}-\Pi_{n,a}^{j-})|> \epsilon_3/a_n \}}\\
&\leq \sup_{|(x,v)|< \delta} a_n^2 n^{-2d_1}  \ \sum_{a=1}^n \mu_{a}(x,v) 2 \mathds{1}_{\{|n^{- d_1} |> \epsilon_3/a_n \}}\\
&\leq \sup_{|(x,v)|< \delta}  \left| \frac{\mu_{\mathbb{A}_n}(x,v)}{n}\right| 2 (a_n n^{1/2-d_1})^2 \mathds{1}_{\{a_n n^{- d_1} > \epsilon_3 \}}
\end{split}
\end{equation}

Now $\{a_n n^{- d_1} > \epsilon_3 \} \xrightarrow{n \to \infty} \emptyset$ and $a_n n^{1/2-d_1} < \infty$ since we assumed $a_n \sqrt{n} = O(n^{d_1})$. 
\newpage
\subsection{Proof of Proposition \ref{prop_stat_dist}}\label{proof_prop_stat_dist}

For simplicity we scale the average opinion by $\bar{v}' := \frac{n \bar{v}}{2}$ so it has its values on the lattice $\mathbb{L}' = \frac{n}{2}\mathbb{L} \subset \mathbb{Z}$.
Following Weidlich and Haag \cite{wei1983} chapter 2.3.1. the stationary distribution of $\bar{v}'$ is recursively given by

\begin{equation}\label{Eq_cases_stat_dist1}
\widehat{\p}_{st}(\bar{v}') = \widehat{\p}_{st}(0) \prod_{y=1}^{\bar{v}'} \frac{\left(\frac{n}{2}-(y-1)\right)\Pi^{1,2}\left(\frac{2(y-1)}{n}\right)}{\left(\frac{n}{2}+y\right)\Pi^{2,1}\left(\frac{2y}{n}\right)} , \ \ \ 1 \leq \bar{v}' \leq n/2
\end{equation}
and
\begin{equation}\label{Eq_cases_stat_dist2}
\widehat{\p}_{st}(\bar{v}') = \widehat{\p}_{st}(0) \prod_{y=1}^{\bar{v}'} \frac{\left(\frac{n}{2}-(y+1)\right)\Pi^{2,1}\left(\frac{2(y+1)}{n}\right)}{\left(\frac{n}{2}-y\right)\Pi^{1,2}\left(\frac{2y}{n}\right)} ,  -n/2 \  \leq \bar{v}' \leq -1.
\end{equation}

Inserting Equations (\ref{state1_example_lux}) and (\ref{state2_example_lux}) into (\ref{Eq_cases_stat_dist1}) yields for $1 \leq \bar{v}' \leq n/2$

\begin{equation}\label{simplify_case1_stat_dist}
\begin{split}
\widehat{\p}_{st}(\bar{v}') &= \widehat{\p}_{st}(0) \prod_{y=1}^{\bar{v}'} \frac{\left(\frac{n}{2}-(y-1)\right) }{\left(\frac{n}{2}+y\right)}\exp \left(\frac{2\gamma}{n} (2y-1)\right)\\
&= \widehat{\p}_{st}(0) \left(\prod_{y=1}^{\bar{v}'} \frac{\left(\frac{n}{2}-(y-1)\right) }{\left(\frac{n}{2}+y\right)}\right) \exp \left(\frac{2\gamma}{n} \sum_{y=1}^{\bar{v}'}  2y-1\right)\\
&= \widehat{\p}_{st}(0) \frac{\left(\frac{n}{2} !\right)^2}{n!} \binom{n}{\frac{n}{2} + \bar{v}'} \exp\left(\frac{2\gamma}{n} (\bar{v}')^2\right),
\end{split}
\end{equation}
where we used that $\sum_{y=1}^{\bar{v}'}  2y-1 = (\bar{v}')^2$ and $\prod_{y=1}^{\bar{v}'} \frac{\left(\frac{n}{2}-(y-1)\right) }{\left(\frac{n}{2}+y\right)} = \frac{\left(\frac{n}{2} !\right)^2}{n!} \binom{n}{\frac{n}{2} + \bar{v}'}$. By symmetry the last representation in Equation (\ref{simplify_case1_stat_dist}) is also true for $-n/2 \  \leq \bar{v}' \leq -1$ and thus after re-scaling we have

\begin{equation}\label{App_P_st}
\p_{st}(\bar{v}) = \widehat{\p}_{st}(\frac{n \bar{v}}{2}) = \p_{st}(0) \frac{\left(\frac{n}{2} !\right)^2}{n!} \binom{n}{\frac{n(1+\bar{v})}{2}} \exp\left(\frac{\gamma n \bar{v}^2}{2}\right), \bar{v} \in \mathbb{L}.
\end{equation}

Next we derive the requirement on $\gamma$ in order that $\p_{st}$ has a local maximum. By symmetry $\p_{st}$ has a local maximum at 0 if the difference at the next higher lattice point is greater 0, i.e.

\begin{equation}
\p_{st}(2/n) - \p_{st}(0) > 0.
\end{equation}

Inserting representation given in Equation (\ref{App_P_st}) yields
\begin{equation}
\begin{split}
&\p_{st}(2/n) - \p_{st}(0) > 0\\
&\Leftrightarrow \frac{\left(\frac{n}{2} !\right)^2}{n!} \binom{n}{\frac{n+2}{2}} \exp\left(\frac{2 \gamma}{n}\right) > 1\\
&\Leftrightarrow \gamma > \frac{n}{2} \ln \left( \frac{n+2}{n} \right) 
\end{split}
\end{equation}

And analogously  $\p_{st}$ has a local minimum for  $\gamma < \frac{n}{2} \ln \left( \frac{n+2}{n} \right) $.

\subsection{Proof of Proposition \ref{DIff_example_lux1}}\label{proof_DIff_example_lux1}

Since $\overline{M}_k = M^2_k - M^1_k \ \forall k \geq 0$, the aggregated state transition for states $s_1 = -1$ and $s_2 = 1$ are given by Equations (\ref{aggregate_minus}), (\ref{aggregate_plus}), (\ref{b_n^i}), (\ref{state1_example_lux}) and (\ref{state2_example_lux}) as

\begin{equation}
\begin{split}
b_n^1(v) &= n^{-1} \sum_{a=1}^n \mu_a(v) \left(\Pi^{1+}(v) - \Pi^{1-}(v) \right)\\
&= v_2 \beta e^{- \gamma (v_2-v_1)} - v_1 \beta e^{ \gamma (v_2-v_1)}\\
&= - \beta \left[ \tanh(\gamma (v_2-v_1)) - (v_2-v_1) \right] \cosh(\gamma (v_2-v_1))
\end{split}
\end{equation}
and analogously
\begin{equation}
\begin{split}
b_n^2(v) = \beta \left[ \tanh(\gamma (v_2-v_1)) - (v_2-v_1) \right] \cosh(\gamma (v_2-v_1)).
\end{split}
\end{equation}

Moreover the transition volume is given by the functions

\begin{equation}
\begin{split}
c_n^i(v) &=n^{-2} \sum_{a=1}^n \mu_a(v) \left(\Pi^{i+}(v) + \Pi^{i-}(v) \right)\\
&= \frac{1}{n} \left( v_2 \beta e^{-\gamma (v_2 - v_1)} + v_1 \beta e^{\gamma (v_2 - v_1) }\right) \\
&= \frac{\beta}{n} \left[ 1 - (v_2-v_1) \tanh(\gamma (v_2-v_1)) \right]\cosh(\gamma (v_2-v_1)) 
\end{split}
\end{equation}

and analogously

\begin{equation}
\begin{split}
c_n^{i,j}(v) = - \frac{\beta}{n} \left[ 1 - (v_2-v_1) \tanh(\gamma (v_2-v_1)) \right]\cosh(\gamma (v_2-v_1)). 
\end{split}
\end{equation}

As $b^1_n$ and $b^2_n$ are independent of $n$ we set $b^1 := b^1_n$ and $b^2 := b^2_n$. Moreover, since $|c_n ^i(v)| = |c_n ^{i,j}(v)| \leq 2 n^{-1}$ we have 
$c_n(v)  \xrightarrow{n \to \infty} 0$.

Now by Theorem \ref{Diff_approx} we get that

\begin{equation}
(V_t^n)_{t \in [0, \infty)} \xrightarrow{\mathcal{L}}  (V_t)_{t \in [0, \infty)} \ in \ D_{[0,1]^2}[0,\infty),
\end{equation}
where $(V_t)_{t \in [0, \infty)}$ is the unique strong solution of
\begin{equation}\label{SDE_ex_2}
\begin{cases}
dV^1_t =  - \beta\left( \tanh(\gamma (V^2_t-V^1_t)) - V^2_t + V^1_t \right)\cosh(\gamma (V^2_t-V^1_t))dt\\
dV^2_t = -dV^1_t
\end{cases} , V_0 = v_0
\end{equation}

which is equivalent to $(\overline{V}_t)_{t \in [0, \infty)}$ being the unique strong solution of
\begin{equation}
d\overline{V}_t =  2 \beta\left( \tanh(\gamma \overline{V}_t) - \overline{V}_t \right)\cosh(\gamma \overline{V}_t) , V_0 = \overline{\theta},
\end{equation}

if we set $\overline{V}_t = V_t^2 - V_t^1$ and $\overline{\theta} = v_0^2 - v_0^1$.


 \newpage
\newpage
\subsection{Proof of Proposition \ref{Diff_ex2}}\label{proof_DIff_ex2}

To apply Theorem \ref{Diff_approx} we calculate the functions $z_n$, $b_n$, $\sigma_n$ and $c_n$ and show their convergence when $n \to \infty$.

Since 
\begin{equation}
\sum^n_{\substack{a \in \mathbb{A}_n, \\ x^a_k \in \{-1,1\}}} e^n_a(P_{k-1},\xi_k) = \sum^n_{\substack{a \in \mathbb{A}_n, \\ x^a_k \in \{-1,1\}}} \widetilde{e}^n_a(P_{k-1},M_{k-1}, \xi_k), \ \forall k \geq 0,
\end{equation}
with
\begin{equation}
\widetilde{e}^n_a(P_{k-1},M_{k-1},\xi_k)=
\begin{cases}
n^{-1/2} w_2(F-P_{k-1}) + \xi_k, & x_k^a = 2\\
n^{-1/2} w_1 \overline{M}_{k-1} + \xi_k, & x^a_k \in \{-1,1\}
\end{cases}
\end{equation}
we have
\begin{equation}
\begin{split}
z_n(x,v) &= n^{-d_2} \sum_{a=1}^n \lambda_a(x,v) \E [e_a^n(x,v,\xi_1)]\\
&= \frac{\overline{\lambda}}{\sqrt{n}} \left( \sum^n_{\substack{a \in \mathbb{A}_n, \\ x^a_0 \in \{-1,1\}}} \E[e^n_a(x,s)] + \sum^n_{\substack{a \in \mathbb{A}_n, \\ x^a_0 =2}} \E[e^n_a(x,s)] \right)\\
&= \frac{\overline{\lambda}}{\sqrt{n}} \left( \sum^n_{\substack{a \in \mathbb{A}_n, \\ x^a_k \in \{-1,1\}}} \E[\widetilde{e}^n_a(x,v,s)] + \sum^n_{\substack{a \in \mathbb{A}_n, \\ x^a_0 =2}} \E[e^n_a(x,s)] \right)\\
&= \frac{\overline{\lambda}}{\sqrt{n}} \left( n(1-\phi_n) n^{-1/2} w_1\bar{v} + n \phi_n n^{-1/2} w_2 (F-x) \right)\\
&=\overline{\lambda} \left( (1-\phi_n)w_1\bar{v} + \phi_n w_2 (F-x) \right)\\
& \xrightarrow{n \to \infty} \overline{\lambda} \left( (1-\phi)w_1\bar{v} + \phi w_2 (F-x) \right) =: z(x,v).
\end{split}
\end{equation}

Moreover

\begin{equation}
\begin{split}
b_n^1(x,v) &= v_2 \Pi^{2,1}_n(x,v) - v_1 \Pi^{1,2}_n(x,v)\\
&= (1-\phi_n)\frac{1+\bar{v}}{2}  \beta e^{-\gamma_1 \widehat{z_n}(x,v) - \gamma_2 \bar{v}} - (1-\phi_n)\frac{1-\bar{v}}{2} \beta e^{\gamma_1 \widehat{z_n}(x,v) - \gamma_2 \bar{v}} \\
&= - (1-\phi_n) \beta \left[ \tanh(\gamma_1 \widehat{z_n}(x,v) - \gamma_2 \bar{v}) - \bar{v} \right] \cosh(\gamma_1 \widehat{z_n}(x,v) - \gamma_2 \bar{v})\\
& \xrightarrow{n \to \infty} - (1-\phi) \beta \left[ \tanh(\gamma_1 z(x,v) - \gamma_2 \bar{v}) - \bar{v} \right] \cosh(\gamma_1 z(x,v) - \gamma_2 \bar{v}),
\end{split}
\end{equation}

since $\phi_n \to \phi$ and $\widehat{z_n}(x,v) \to z(x,v)$ for $n \to \infty$.

\begin{equation}
\begin{split}
b_n^2(x,v) &= v_1 \Pi^{1,2}_n(x,v) - v_2 \Pi^{2,1}_n(x,v)\\
&= -b_n^1(x,v)\\
& \xrightarrow{n \to \infty} (1-\phi) \beta \left[ \tanh(\gamma_1 z(x,v) - \gamma_2 \bar{v}) - \bar{v} \right] \cosh(\gamma_1 z(x,v) - \gamma_2 \bar{v}).
\end{split}
\end{equation}

Furthermore,

\begin{equation}
\begin{split}
&\sigma_n(x,v)^2= \frac{1}{n} \sum_{a=1}^n \bar{\lambda}\E [e_a^n(x,s)^2]\\
&= \frac{\bar{\lambda}}{n}  \left( \sum^n_{\substack{a \in \mathbb{A}_n, \\ x^a_0 \in \{-1,1\}}} \E [e_a^n(x,s)^2] + \sum^n_{\substack{a \in \mathbb{A}_n, \\ x^a_0 = 2}} \E [e_a^n(x,s)^2] \right) \\
&= \frac{\bar{\lambda}}{n} [ (n -k_n) n^{-1} w_1^2 + k_n n^{-1} w_2^2 (F-x)^2 + n\Var[s] ]\\
&= \frac{\bar{\lambda}}{n} [ (1 -\phi_n) w_1^2 + \phi_n w_2^2 (F-x)^2 + n\Var[s]]\\
&\xrightarrow{n \to \infty} \bar{\lambda} \Var[s].
\end{split}
\end{equation}
At last, we have
\begin{equation}
\begin{split}
c_n^i(x,v)^2 &= n^{-2} \sum_{a=1}^n ( \Pi^{i+}_n(x,v) + \Pi^{i-}_n(x,v))\\
&= n^{-1} ( \Pi^{i+}_n(x,v) + \Pi^{i-}_n(x,v))\\
&\leq 2/n \xrightarrow{n \to \infty} 0
\end{split}
\end{equation}
as well as
\begin{equation}
\begin{split}
c_n^{i,j}(x,v)^2 &= -n^{-2} \sum_{a=1}^n ( v_i\Pi^{i,j}_n(x,v) + v_j\Pi^{j,i}_n(x,v))\\
&= -n^{-1} ( v_i\Pi^{i,j}_n(x,v) + v_j\Pi^{j,i}_n(x,v))\\
&\leq -2/n \xrightarrow{n \to \infty} 0.
\end{split}
\end{equation}

Now, by Theorem \ref{Diff_approx} we have

\begin{equation}
(X_t^n,V_t^n)_{t \in [0, \infty)} \xrightarrow{\mathcal{L}}  (X_t,V_t)_{t \in [0, \infty)} \ in \ D_{\mathbb{R} \times [0,1]^3}[0,\infty),
\end{equation}
where $(X_t,V_t)_{t \in [0, \infty)}$ is the unique strong solution of
%

\begin{equation}
\begin{cases}
dX_t = \alpha \overline{\lambda}[\phi w_2 (F-X_t) + (1-\phi) w_1 \overline{V}_t] dt + \sqrt{\bar{\lambda} \Var[s]} dB_t , &X_0 = \eta\\
d\overline{V}_t = 2 \beta \left[ \tanh ( \gamma_1 dX_t + \gamma_2 \overline{V}_t) - \overline{V}_t \right] \cosh(\gamma_1 dX_t + \gamma_2 \overline{V}_t)dt, & \overline{V}_0 = \overline{\theta}\\
dV_t^3 = 0, & V^3_0 = \phi,
\end{cases}
\end{equation}

when setting $\overline{V}_t = \frac{(V_t^2-V_t^1)}{1-\phi}$.
\newpage
\thispagestyle{plain}
\nocite{*}
\bibliography{Literature}
\bibliographystyle{abbrv}







\end{document}